\documentclass[11pt]{article}

\usepackage{shapepar,bm}
\usepackage[dvips,arrow,matrix,ps,color,line]{xy}
\usepackage{theorem,amsfonts,amssymb,amscd}
\usepackage{geometry}

\usepackage{graphicx}
\usepackage{graphics}
\usepackage{amssymb}
\usepackage{epstopdf}

\usepackage{amssymb, graphics,hyperref}

\hypersetup{
    colorlinks = true,
    linkcolor = blue,
    anchorcolor = red,
   citecolor = blue,
    filecolor = red,
    pagecolor = red,
    urlcolor = blue}

\newcommand{\mathsym}[1]{{}}

\usepackage[usenames]{color}
\definecolor{MyLightMagenta}{cmyk}{0.1,0.8,0,0.1}
\definecolor{MyDarkBlue}{rgb}{0.1,0,0.3}

\makeindex

\def\Lcal{\mathcal L}

\def\NN{\mathbb N}

\def\bfe{{\mathbf e}}
\def\bfz{{\mathbf z}}
\def\bfa{{\mathbf a}}

\def\bfc{{\mathbf c}}
\def\bfh{{\mathbf h}}
\def\bfv{{\mathbf v}}
\def\Qcal{\mathcal Q}

\def\Bcal{{\mathcal B}}
\def\Ccal{{\mathcal C}}

\def\Ucal{\mathcal U}
\def\Vcal{\mathcal V}
\def\Mcal{\mathcal M}
\def\Ncal{\mathcal N}
\def\Ocal{\mathcal O}
\def\Gcal{\mathcal G}
\def\Uab{U_\alpha\cap U_\beta}
\def\ab{{\alpha\beta}}
\def\bfx{{\mathbf x}}

\def\bff{{\mathbf f}}

\def\mod{{\mathrm{mod}}}
\def\Sym{{\mathrm{Sym}}}

\def\ZZ{\mathbb Z}
\def\EE{\mathbb E}

\def\CC{\mathbb C}

\def\bff{{\bf f}}

\def\wlbund{L^{\otimes r+1}\otimes K^{\otimes{r(r+1)\over 2}}}

\def\KK{\mathbb K}
\def\gf{\pi:\TS\lra S}
\def\RR{\mathbb R}
\def\Scal{{\mathcal S}}
\def\Ecal{{\mathcal E}}
\def\WW{\mathbb W}
\def\QQ{\mathbb Q}
\def\PP{\mathbb P}
\def\cocoa{{\hbox{\rm C\kern-.13em o\kern-.07em C\kern-.13em o\kern-.15em A}}}

\def\Dcal{\mathcal D}

\def\Acal{{\mathcal A}}
\def\Fcal{{\mathcal F}}

\def\Grd{\Gamma(\rho_{r,d})}
\def\Gtrdb{\Gamma_{\tt triv}(\rho_{r,d})}

\def\Poly{{\rm Poly}}
\def\Res{{\rm Res}}

\def\F{\mathcal{F}}

\def\Ical{\mathcal I}

\def\rk{{\mathrm{rk}}}

\def\TS{{\frak X}}

\def\bfu{{\bf u}}

\def\partoformu{{\mu^{\lambda_r}\wedge\mu^{1+\lambda_{r-1}}\wedge\ldots\wedge\mu^{r+\lambda_0}}}
\def\partorformep{{\ep^{\lambda_r}\wedge\ep^{1+\lambda_{r-1}}\wedge\ldots\wedge\ep^{r+\lambda_0}}}

\def\blamb{{\bm \lambda}}
\def\bmu{{\bm \mu}}

\def\Pcal{{\mathcal P}}

\def\ep{{\epsilon}}

\def\wt{{\rm wt}}
\def\ev{{\rm ev}}

\def\w2M{\bigwedge^2M}

\def\w{\wedge }
\def\bw{\bigwedge }

\protect
\protect
\protect
\protect

\protect
\protect

\def\sra{\rightarrow}
\def\lra{\longrightarrow}

\def\proof{\noindent{\bf Proof.}\,\,}
\def\qed{{\hfill\vrule height4pt width4pt depth0pt}\medskip}
\def\be{\begin{equation}}
\def\ee{\end{equation}}
\def\bclm{\begin{claim}}
\def\eclm{\end{claim}}
\def\beqn{\begin{eqnarray}}
\def\eeqn{\end{eqnarray}}
\def\beqn*{\begin{eqnarray*}}
\def\eeqn*{\end{eqnarray*}}

\theoremstyle{change}
\theorembodyfont{\rmfamily}

\newtheorem{claim}{}[section]

\global
\global

\def\no@breaks#1{{\def\\{ \ignorespaces}#1}}    

  \makeatletter
\def\cleardoublepage{\clearpage\if@twoside \ifodd\c@page\else
\hbox{} \thispagestyle{empty}
\newpage
\if@twocolumn\hbox{}\newpage\fi\fi\fi} \makeatother

\usepackage{eso-pic,graphicx}
\makeatletter
\newcommand\BackgroundPicture[2]{%
  \setlength{\unitlength}{1pt}%
  default \put(0,\strip@pt\paperheight){%
  \parbox[t][\paperheight]{\paperwidth}{%
    \vfill
     \centering \includegraphics[angle=#2, width=15cm, height=15cm,  bb=0 0 150 150]{#1}
    \vfill
}}} %
\makeatother

\usepackage{palatino,url}  

\title{\Large  On Generalized Wro\'nskians{}\thanks{\noindent Work partially sponsored by PRIN ``Geometria
sulle Variet\`a Algebriche" (Coordinatore A.~Verra), Politecnico di Torino. The second author was
sponsored by an INDAM-GNSAGA grant (2009) for Visiting Professors at the Politecnico di Torino.}\\ }

\author{\large L.~Gatto\,\,\,\,\,\,\,\,\,\,\,\qquad\quad\qquad\qquad\,\,\,\,\,\,\,\,\,\,\, I.~Scherbak \\{}\\
{ Dipartimento di Matematica\qquad\qquad
  School of Mathematical Sciences }\\
{ Politecnico di Torino\,\,\,\,\,\,\,\,\,\,\,\,\,\,\,\, \qquad\qquad\qquad  Tel Aviv
University\,\,\,\,\, }}

\date{}                                           

\begin{document}

\maketitle

\noindent { \abstract{\noindent The   Wro\'nski determinant ({\em Wro\'nskian}),  usually introduced in
standard courses in Ordinary Differential Equations (ODE), is a very useful tool  in algebraic geometry
to detect ramification loci of linear systems. The present survey aims to describe some
"materializations" of the Wro\'nskian and of its close relatives, {\it the generalized Wro\'nskians}, in
algebraic geometry. Emphasis will be put  on the relationships between Schubert Calculus and ODE.}}



\tableofcontents

\section*{Introduction}
 \addcontentsline{toc}{section}{Introduction}

Let $\bff:=(f_0,f_1,\ldots, f_r)$ be an $(r+1)$-tuple of holomorphic functions in one complex variable.
The {\em Wro\'nskian} of $\bff$ is the holomorphic function $W(\bff)$ obtained by taking the determinant
of the {\em Wro\'nski matrix} whose entries of the $j$-th-row, $0\leq j\leq r$, are the $j$-th
derivatives of $(f_0,f_1,\ldots, f_r)$. The first appearance of Wro\'nskians dates back to 1812,
introduced  by J.~M.~Hoene--Wro\'nski (1776--1853) in the treatise~\cite{HWr} -- see
also~\cite{pragacz2}. The ubiquity of the  Wro\'nskian  in nearly all the branches of mathematics, from
analysis  to algebraic geometry, from number theory to combinatorics, up to the theory of infinite
dimensional dynamical systems,  is definitely surprising if compared with its   elementary  definition.
The present survey aims to outline links between some different Wro\'nskian materializations  to make
evident their common root. The emphasis will be put on  the mutual relationships among linear Ordinary
Differential Equations (ODEs), the theory of ramification loci of linear systems (e.g. Weierstrass
points on curves)  and the intersection theory of complex Grassmann varieties, ruled by the famous  {\em
Calculus}  \cite{schubert} elaborated in 1886~by H.~C.~H.~Schubert  (1848--1911), to which the Italians
M.~Pieri (1860--1913) and G.~Z.~Giambelli (1879--1953) contributed too -- see \cite{giambelli, pieri}.

The notion of Wro\'nskian belongs to mathematicians' common background   because of  its most popular
application, which provides  a method (sketched in Section~\ref{walodes})  to   find a particular
solution of a non-homogeneous linear ODE. It  relies on the following  key property of the Wro\'nskian
of a  fundamental system of solutions of a linear homogeneous ODE: {\em the derivative of the
Wro\'nskian is proportional to the Wro\'nskian itself}, whose proof  is due to J.~Liouville (1809--1882)
and N.~H.~Abel (1802--1829). This apparently innocuous property should be considered  as the first
historical appearance  of Schubert Calculus. To see it, one must embed the Wro\'nski determinant into a
full family of {\em generalized Wro\'nskians}, already used in 1939 by F.~H.~Schmidt \cite{Schmidt} to
study Weierstrass points and, in recent times and with the same motivation,    by C.~Towse
in~\cite{towse}. For a sample  of applications to number theory see also~\cite{lacunary} and
~\cite{MMO}.

If $\blamb=(\lambda_0\geq\lambda_1\geq \ldots\geq \lambda_r)$ is a {\em partition}, the {\em generalized
Wro\'nskian} $W_\blamb(\bff)$ is the determinant of the matrix whose $j$-th row, for $0\leq j\leq r$, is
the row of the derivatives of order $j+\lambda_{r-j}$ of $(f_0,f_1,\ldots, f_r)$. Clearly
$W(\bff)=W_0(\bff)$, where $0$ stands for the {\em null partition} $(0,\ldots,0)$. The derivative of
$W(\bff)$, appeared in the proof of  Liouville's--Abel's theorem, is the first example of a generalized
Wro\'nskian, $W_{(1)}(\bff)$, corresponding to the partition $(1,0,\ldots, 0)$. The bridge to Schubert
Calculus is our generalization of Liouville's and Abel's theorem (see \cite{GatScherb1}): {\em
Giambelli's formula for generalized Wro\'nskians holds}. More precisely, if $\bff$ is a fundamental
system of solutions of a linear ODE {\em with constant coefficients}, then  $W_\blamb(\bff)$ is
proportional to the usual Wro\'nski determinant, $W_\blamb(\bff)=\Delta_\blamb(\bar\bfh) W_0(\bff)$,
where $\Delta_\blamb(\bar\bfh)$ is the Schur polynomial associated to a sequence
$\bar\bfh=(h_0,h_1,\ldots)$ of explicit  polynomial expressions in the coefficients of the given ODE and
to the partition $\blamb$ -- see Section~\ref{finalsec}. If the characteristic polynomial of the linear
differential equation splits into the product of distinct linear factors, then  $h_j$ is nothing else
than the $j$-th complete symmetric polynomial in its roots.

\medskip
Let us now change the landscape for a while.  Take a smooth complex projective curve $C$ of genus $g\geq
0$ and an isomorphism class $L\in Pic^d(C)$ of line bundles of degree $d$ on $C$. A  $g^r_d$ on $C$ is a
pair $(V,L)$, where $V$ is a point of the  Grassmann variety $G(r+1, H^0(L))$ parameterizing of
$(r+1)$-dimensional vector subspaces of the global holomorphic  sections of $L$. If
$\bfv=(v_0,v_1,\ldots, v_r)$ is a basis of $V,$ the {\em Wro\'nskian} $W(\bfv)$  is a holomorphic
section of the bundle $\Lcal_{g,r,d}:=\wlbund$ -- see Section~\ref{wrseclbclass}. It can be constructed
by gluing together local Wro\'nskians $W(\bff)$, where $\bff=(f_0,f_1,\ldots,f_r)$  is an $(r+1)$-tuple
of holomorphic functions representing the basis $\bfv$ in some open set of $C$ that trivializes $L$.  As
changing the basis of $V$ amounts to multiply $W(\bfv)$ by a non-zero complex number, one obtains a well
defined point $W(V):=W(\bfv)\,(\mod\,\CC^*)$ in $\PP H^0(\Lcal_{g,r,d})$ called  the {\em Wro\'nskian}
of $V$. The {\em Wro\'nski map} $G(r+1, H^0(L))\sra \PP H^0(\Lcal_{g,r,d})$ mapping $V$ to $W(V)$  is a
holomorphic map; two extremal cases show that, in general, it  is neither injective nor surjective.
Indeed, if   $C$ is  hyperelliptic and $L\in Pic^2(C)$ is the line bundle defining its unique $g^1_2$,
then $G(2, H^0(L))$ is just a point, and the Wro\'nski map to $\PP H^0(\Lcal_{g,1,2})$ is trivially
injective and not surjective. On the other hand,
 if $C=\PP^1$ and $L=O_{\PP^1}(d)$, then the Wro\'nski map $G(r+1, H^0(O_{\PP^1}(d)))\sra \PP H^0(\Lcal_{0,r,d})$
 is a finite surjective morphism whose degree is equal to the Pl\"ucker degree of the Grassmannian
 $G(r+1,d+1)$, thence in this case the Wro\'nski map is not injective, cf. \cite{EHcuspidal}.

\medskip
 The problem of determining the pre-image of an element of  $\PP H^0(\Lcal_{0,r,d})$ through the Wro\'nski map
 defined on $G(r+1, H^0(O_{\PP^1}(d)))$ leads to an intriguing  mixing of Geometry,  Analysis and Representation
 Theory. It turns out that certain {\em non-degenerate elements} of $G(r+1, H^0(O_{\PP^1}(d)))$,
 defined through suitable {\em intermediate Wro\'nskians},  correspond to the so-called Bethe vectors appeared
 in  representation theory
 of the Lie algebra ${\frak s l}_{r+1}(\CC)$. The correspondence goes through critical points of a
 remarkable rational function related to Knizhnikov-Zamolodchikov equation on correlation functions
 of the conformal field theory,  ~\cite{MV,Scherb1,Scherb2,Scherb3}.
Interestingly, the critical points of the mentioned rational function in the case $r=1$ were examined in
the XIX century, in works of Heine and Stieltjes on second order Fuchsian differential equations having
a polynomial solution of a prescribed degree. Schubert Calculus on Grassmannians has been introduced
even before. However, the relationship between these items - in the case $r=1$ - was conceived a decade
ago in~\cite{Scherb,Scherb3}.

In the real framework, the relationship between Wro\'nskians, Schubert Calculus and rational curves was
discovered and studied by L.~Goldberg, A.~Eremenko\& A.~Gabrielov,  V.~Karlhamov \& F.~Sottile, and
others -- see \cite{ Gold, EreGab1, EreGab2, KhSo} and references therein.

 More links between linear differential equations, projective curves and Schubert varieties
 appeared in a local context in the investigations of M.~Kazarian on singularities of the boundary
of fundamental systems of solutions of linear differential equations, \cite{K}.

\medskip
Here, we take another point of view.   A.~Nigro proposes  to extend the notion of ramification locus of
a linear system on a curve to that of ramification locus of a holomorphic section of a Grassmann bundle
\cite{NigroTesi}. The construction was motivated by the following observation (see also~\cite{CuGaNi}):
Let $\Gamma_{\tt triv}(\rho_{r,d})$ be the set of all the sections $\gamma:C\sra G(r+1, J^dL)$  such
that the pull back of the tautological bundle $\Scal_r$ over $G(r+1, J^dL)$ is trivial.  Then each
$g^r_d:=(V,L)$ induces a holomorphic section $\gamma_V\in \Gamma_{\tt triv}(\rho_{r,d})$, via the bundle
monomorphism $C\times V\sra J^dL$ (cf. Section~\ref{sec53}). The point is that the space $\Gamma_{\tt
triv}(\rho_{r,d})$ is larger than the space of linear systems, and so the theory becomes richer. A
distinguished subvariety  indwells in $G(r+1, J^dL)$,  called {\em Wro\'nski subvariety}
in~\cite{NigroTesi}. It is a Cartier divisor which occurs as the zero locus of a certain {\em Wro\'nski
section} $\WW$. The Wro\'nskian of any section $\gamma\in\Gamma_{\tt triv}(\rho_{r,d})$ is defined to be
$W_0(\gamma):=\gamma^*\WW \,(\mod\,\CC^*)$;   if  $\gamma=\gamma_V$ for some $V\in G(r+1, H^0(L))$, it
coincides with the usual Wro\'nskian of  $V$ -- see Section~\ref{wsgb}. In particular, if $\Mcal$ is  a
line bundle defining the unique $g^1_2$ over a hyperelliptic curve of genus $g\geq 2$, the extended
Wro\'nski map $\Gamma_{\tt triv}(\rho_{1,2})\sra \PP H^0(\Mcal^{\otimes 2}\otimes K)$ is dominant
(see~\cite{CuGaNi}), its behavior is closer   to  the surjectivity of the Wro\'nski map defined on the
space of  $g^r_d$s on $\PP^1$. The latter, in this case,  coincides with $\Gamma_{\tt triv}(\rho_{r,d})$
modulo identification of $V$ with $\gamma_V$.

In general, the  construction  works  as follows.  Let $\varrho:F\sra X$ be a vector bundle of rank
$d+1$  and  $\varrho_{r,d}:G(r+1, F)\sra X$  be the Grassmann bundle of $(r+1)$-dimensional subspaces of
fibers of $\varrho$. Consider $0\sra \Scal_r\sra\varrho_{r,d}^*F\sra \Qcal_r\sra 0$, the universal exact
sequence over $G$, and denote by $\Delta_\blamb(c_t(\Qcal_r-\varrho_{r,d}^*F))$  the Schur polynomial,
associated to the partition $\blamb$,  in the coefficients of the Chern polynomial of
$\Qcal_r-\varrho_{r,d}^*F$. As is well known (see e.g.~\cite[Ch.~14]{Fu1}), the Chow group $A^*(G)$ of
cycles modulo rational equivalence is a free $A^*(X)$-module generated by
$\Bcal:=\{\Delta_\blamb(c_t(\Qcal_r-\varrho_{r,d}^*F))\cap[G]\,|\, \blamb\in\Pcal^{(r+1)\times
(d-r)}\}$, where $\Pcal^{(r+1)\times (d-r)}$ denotes  the set of the partitions $\blamb$ such that
$\lambda_0\leq d-r$,  and $\cdot\cap[G]$ denotes the {\em cap} product with the fundamental class of
$G$. Let $F_\bullet:=(F_i)_{d\geq i\geq  0}$ be a filtration of $F_\bullet$ by quotient bundles, such
that $F_i$ has rank $i$. Schubert varieties $\{\Omega_\blamb(F_\bullet)\,|\, \blamb\in\Pcal^{(r+1)\times
(d-r)}\}$ associated to $F_\bullet$ (the definition is in Section~\ref{svf44}) play the role of  {\em
generalized Wro\'nski subvarieties}.
 In particular $\Omega_{(1)}(F_\bullet)$ is what in~\cite{NigroTesi} was called the $F_\bullet$-{\em Wro\'nski
subvariety} of $G$. It is a Cartier divisor, that is the zero locus of a section $\WW$ of the bundle
$\bw^{r+1}\varrho_{r,d}^*F_{r+1}\otimes \bw^{r+1}\Scal_r^\vee$ over $G$. We say that  $\WW$ is the
$F_\bullet$-{\em Wro\'nskian}. If $\gamma:X\sra G$ is a holomorphic section, its Wro\'nskian is, by
definition, $W_0(\gamma):=\gamma^*\WW\in H^0(\bw^{r+1}F\otimes\bw^{r+1}\gamma^*\Scal_r^\vee)$. Its class
in $A_*(X)$ is nothing else than $\gamma^*[\Omega_{(1)}(F_\bullet)]\cap [G]$. The {\em generalized
Wro\'nski class}   of $\gamma$ in $A_*(X)$ is $\gamma^*[\Omega_\blamb(F_\bullet)]\cap[X]$, which is the
class of $\gamma^{-1}(\Omega_\blamb(F_\bullet))$, provided that the codimension of the locus coincides
with the expected codimension $|\blamb|:=\lambda_0+\ldots +\lambda_r$. Recall that
$[\Omega_\blamb(F_\bullet)]$ can be easily computed as an explicit linear combination of the elements of
the  basis $\Bcal$ above, for instance by the recipe indicated in Section~\ref{wrsecbngen}, especially
Theorem~\ref{kmpflak}.

Let now $\varepsilon_i:=c_i(\Scal_r)\in A^*(G)$ be the Chern classes of the tautological bundle
$\Scal_r\sra G$. Consider a basis $\bfv:=(v_0,v_1,\ldots, v_r)$ of solutions of the differential
equation
\be y^{(r+1)}-\varepsilon_1y^{(r)}+\ldots+(-1)^{r+1}\varepsilon_{r+1}y=0,\label{eq:schubdifeq}
\ee
taken in the algebra $(A^*(G)\otimes\QQ)[[t]]$ of formal power series in an indeterminate $t$ with
coefficients in the Chow ring of $G$ with rational coefficients. In Section~\ref{fin813} we show that,
for each partition $\blamb\in\Pcal^{(d-r)\times (r+1)}$,
\[
\Delta_\blamb(c_t(\Qcal_r-\rho_{r,d}^*F))={W_\blamb(\bfv)\over W_0(\bfv)},
\]
i.e. each element of the $A^*(X)$-basis of the Chow ring of $G$  is the quotient of generalized
Wro\'nskians associated to a fundamental system of solutions of an ordinary linear ODE with constant
coefficients taken in $A^*(G)$. This will be a  consequence of {\em Giambelli's formula} for generalized
Wro\'nskians, proven in~\cite{GatScherb1}, which so provides another clue of the ubiquity of
Wro\'nskians in mathematics.

\medskip The survey was written with an eye on a wide range of readers, not necessarily experts in
algebraic geometry. We thank the referees for substantial efforts to improve the presentation.

\section{Wro\'nskians,  in General}

\claim{}\label{s11} In the next two  sections let  $\KK$ be  either the real field  $\RR$ or the complex
field $\CC$  together with their usual euclidean topologies. If $U\subseteq \KK$ is an open connected
subset of $\KK$, we shall write $\Ocal(U)$ for the $\KK$-algebra of {\em regular} $\KK$-valued functions
defined over $U$: here {\em regular} means either  $C^\infty$ differentiable  if $\KK=\RR$ or complex
holomorphic if $\KK=\CC$. Let
\be
\bfv:=(v_0, v_1, \ldots, v_r)\in \Ocal(U)^{r+1}.\label{eq:functions}
\ee
If $t$ is  a local  parameter on $U$, we denote by  $D: \Ocal(U)\sra \Ocal(U)$ the  usual derivation
$d/dt$.  The {\em Wro\'nski matrix} associated to the $(r+1)$-tuple~(\ref{eq:functions}) is the matrix
valued regular function:
 \[
 WM(\bfv):=\pmatrix{\bfv\cr
 D\bfv\cr\vdots\cr D^r\bfv}=\pmatrix{v_0&v_1&\ldots&v_r\cr
 Dv_0&Dv_1&\ldots&Dv_r\cr\vdots&\vdots&\ddots&\vdots\cr  D^rv_0&D^rv_1&\ldots&D^rv_r}.
 \]
The determinant $ W_0(\bfv):=\det (WM(\bfv))$ is the {\em Wro\'nskian }of $\bfv:=(v_0,v_1,\ldots, v_r)$.
It will be often written in the form:
 \be
 W_0(\bfv):=\bfv\w D\bfv\w \ldots\w D^r\bfv.\label{eq:wrrpuform}
 \ee

In this paper, however, we want to see Wro\'nskians as a part of a full family of natural functions
generalizing them. They will be called, following the few pieces of literature where they have already
appeared (\cite{lacunary}, \cite{towse}) {\em generalized Wro\'nskians}.

\claim{\bf Generalized Wro\'nskians.}\label{genwr} Let $r\geq 0$ be an integer. A partition $\blamb$ of
length at most $r+1$ is an $(r+1)$-tuple of non-negative integers in the non-increasing order:
\be \blamb:\quad
\lambda_0\geq\lambda_1\geq \ldots\geq \lambda_r\geq 0.\label{eq:lambpart}
\ee
The {\em weight} of $\blamb=(\lambda_0,\lambda_1,\ldots, \lambda_r)$ is
$|\blamb|:=\sum_{i=0}^r\lambda_j$, that is $\blamb$ is a partition of the integer $|\blamb|$. In this
paper we consider only partitions of length $r+1$. To each partition one may associate a Young--Ferrers
diagram, an array of left justified rows, with $\lambda_0$ boxes in the first row, $\lambda_1$ boxes in
the second row, \ldots, and $\lambda_r$ boxes in the $(r+1)$-th row. We denote by $\Pcal^{(r+1)\times
(d-r)}$ the set of all partitions whose Young diagram is contained in the $(r+1)\times (d-r)$ rectangle,
i.e. the set of all partitions $\blamb$ such that
\[
d-r\geq \lambda_0\geq \lambda_1\geq\ldots\geq \lambda_r\geq 0.
\]
If the last $r- h$ entries of $\blamb\in\Pcal^{(r+1)\times (d-r)}$ are zeros, then we write simply
$\blamb=(\lambda_0,\lambda_1,\ldots,\lambda_h)$, omitting the last zero parts. For more on partitions
see~\cite{MacDonald}.

\bclm{\bf Definition.} {\em Let $\bfv$ as in~(\ref{eq:functions}) and $\blamb$ as
in~(\ref{eq:lambpart}).  The {\em generalized Wro\'nski matrix} associated to $\bfv$ and to  the
partition $\blamb$ is, by definition,
\[
WM_\blamb(\bfv):=\pmatrix{D^{\lambda_r}\bfv\cr D^{1+\lambda_{r-1}}\bfv\cr \vdots\cr
D^{r+\lambda_{0}}\bfv}:=\pmatrix{D^{\lambda_r}v_0& D^{\lambda_r}v_1&\ldots&D^{\lambda_r}v_r\cr
 D^{1+\lambda_{r-1}}v_0& D^{1+\lambda_{r-1}}v_1&\ldots&D^{1+\lambda_{r-1}}v_r\cr\vdots&\vdots&\ddots&\vdots\cr
D^{r+\lambda_0}v_0& D^{r+\lambda_0}v_1&\ldots&D^{r+\lambda_0}v_r}\,.
\]
The $\blamb$-{\em generalized Wro\'nskian} is the determinant of the generalized Wro\'nski matrix:
\[
W_\blamb(\bfv):=\det WM_\blamb(\bfv)\,.
\]}
\eclm
Coherently with~(\ref{eq:wrrpuform}) we shall write the $\blamb$-generalized Wro\'nskian in the form:
\be W_\blamb(\bfv):=D^{\lambda_r}\bfv\w D^{1+\lambda_{r-1}}\bfv\w \ldots\w
D^{r+\lambda_0}\bfv.\label{eq:genwrrpuform} \ee The usual Wro\'nskian corresponds to the partition of
$0$, that is  $W(\bfv)\equiv W_0(\bfv)$.

\claim{\bf Remark.} Notation~(\ref{eq:wrrpuform}) and~(\ref{eq:genwrrpuform}) is convenient because the
derivative of any generalized Wro\'nskian can be computed via  Leibniz's rule with respect to the
product ``$\wedge$'':
\begin{eqnarray*}
D(W_\blamb(\bfv))&=&D(D^{\lambda_r}\bfv\w D^{1+\lambda_{r-1}}\bfv\w \ldots\w
D^{r+\lambda_0}\bfv)=\\&=&\sum_{\scriptsize{\matrix{i_0+i_1+\ldots+i_r=1\cr i_j\geq 0}}}
D^{i_0+\lambda_r}\bfv\w D^{1+i_1+\lambda_{r-1}}\bfv\w \ldots\w D^{r+i_r+\lambda_0}\bfv.
\end{eqnarray*}
A simple induction shows that {\em any} derivative of  $W_\blamb(\bfv)$ is a $\ZZ$-linear combination of
generalized Wro\'nskians.  Recall, as  in Section~\ref{genwr},  that partitions can be described via
Young--Ferrers diagrams, and that a standard {\em Young tableau} is a numbering  of the boxes of the
Young--Ferrers diagram of $\blamb$ with integers $1,\ldots, |\blamb|$ arranged in an increasing order in
each column and each row \cite{Fu2}.

The following observation has convinced us that the Schubert Calculus can be recast in terms of
Wro\'nskians, see Section~\ref{finalsec}.

\bclm{\bf Theorem.} {\em We have
\[
D^{h}W(\bfv)\,=\, \sum_{|\blamb|=h} c_\blamb W_\blamb(\bfv)\,,
\]
where $c_\blamb$ is the number of the standard Young tableaux of the Young--Ferrers diagram $\blamb$.
\qed}
\eclm
The coefficients $c_\blamb$'s and their interpretation in terms of Schubert Calculus are very well
known; in particular, they can be calculated by the {\em hook formula}:
\[
c_\blamb={|\lambda|! \over k_1\cdot\,\ldots\, \cdot k_{|\lambda|}}\,,
\]
where the $k_j$'s, $1\leq j\leq |\blamb|$, are the {\it hook lengths} of the boxes of $\blamb$, see
\cite[p.~53]{Fu2}.

\section{Wro\'nskians and Linear ODEs}\label{walodes}

Wro\'nskians  are usually introduced  when dealing with linear {\em Ordinary Differential Equations}
(ODEs).

\claim{} We use notation of Section~\ref{s11}. For $\bfa(t)=(a_1(t),\ldots, a_{r+1}(t))\in
\Ocal(U)^{r+1}$ and $f\in\Ocal(U)$, consider the linear ODE
\be D^{r+1}x-a_1(t)D^rx+\ldots+(-1)^{r+1}a_{r+1}(t)x=f\label{eq:typODE} \ee
and  the corresponding linear differential operator $P_\bfa(D)\in End_\KK(\Ocal(U))$,
\be
P_\bfa(D):=D^{r+1}-a_1(t)D^{r}+\ldots+(-1)^{r+1}a_{r+1}(t).\label{eq:DifOp}
\ee
The set of solutions, $\Scal_{f,\bfa}$, of~(\ref{eq:typODE})  is an affine space modelled over
$\KK^{r+1}$\,: \ if $x_p$ is a {\em particular solution} , then
\[
\Scal_{f,\bfa}=x_p+\ker P_\bfa(D)\,.
\]
The celebrated Cauchy theorem ensures  that given a column $\bfc=(c_j)_{0\leq j\leq r}\in\KK^{r+1}$,
 there exists a unique element $x_\bfc\in \ker P_\bfa(D)$ such that  $D^jf(0)=c_j$, for all $0\leq j\leq r$.
 Assume now that $\bfv$ as in~(\ref{eq:functions})
is a basis of $\ker P_\bfa(D)$. A particular solution of~(\ref{eq:typODE}) can be found through the
method of {\em variation of arbitrary constants}. Assume that
\[
\bfc=\bfc(t)=\pmatrix{c_0(t)\cr c_1(t)\cr\vdots\cr c_r(t)}\,\in\,\Ocal(U)^{r+1}\,,
\] and look for a solution of~(\ref{eq:typODE}) of the form
\[
x_p:=(\bfv\cdot \bfc)(t)=\bfv(t)\cdot \bfc(t)=\sum_{i=0}^rc_i(t)v_i(t)\,,
\]
where `$\cdot$'' stands for the usual row-by-column product. The condition that $D^j\bfv\cdot D\bfc=0$
for all $0\leq j\leq r$  means that $D^jx_p=D^j\bfv\cdot\bfc$ for all $0\leq j\leq r$ and
$D^{r+1}x_p=D^{r+1}\bfv\cdot \bfc+D^r\bfv \cdot D\bfc$. The equality
\[
P_{\mathbf a}(D)x_p=f
\]
implies, by substitution, the equation
\[
D^r\bfv\cdot D\bfc=f.
\]
The  unknown functions $\bfc=\bfc(t)$ must  then satisfy the differential equations:
\[
WM(\bfv)\pmatrix{Dc_0\cr Dc_1\cr \vdots\cr Dc_r}=\pmatrix{0\cr 0\cr \vdots\cr f}.
\]
The key remark is that the {\em Wro\'nski matrix is invertible in $\Ocal(U)$}. Thus we get a system of
first order ODEs,
\[
D\bfc=(WM(\bfv))^{-1}\cdot\pmatrix{0\cr 0\cr \vdots\cr  f}\,,
\]
which can be solved by usual methods.

 To show the
invertibility, one usually shows that if the Wro\'nski matrix does not vanish at some point of $U$, then
it does vanish nowhere on $U$ (recall that $U$ is a connected open set). Assume  $W_0(\bfv)(P)\neq 0$
for some $P\in U$. Let us choose a local parameter $t$ on $U$ which is $0$ at $P$, identifying the open
set $U$ with a connected neighborhood of the origin. Computing the derivative of the Wro\'nskian one
discovers the celebrated

\bclm{\bf Liouville's Theorem (\cite[p.~195, \S 27.6]{arnold}).}\label{liouvthm} {\em The Wro\'nskian
$W=W_0(\bfv)$ satisfies the differential equation: \be DW=a_1W.\label{eq:liouvthm} \ee } \eclm

The proof of  Theorem~\ref{liouvthm} is as follows. By defining $D\bfv$ as the row whose entries are the
derivatives of the entries of $\bfv$, one notices that
\[
P_\bfa(D)\bfv=(P_\bfa(D)v_0,P_\bfa(D)v_1,\ldots, P_\bfa(D)v_r)={\mathbf 0}\,.
\]
Thence $D^{r+1}\bfv=a_1(t)D^r\bfv-a_2(t)D^{r-1}\bfv+\ldots+(-1)^{r}a_{r+1}(t)\bfv$ and one gets
\begin{eqnarray*}
DW_0(\bfv)&=&D(\bfv\w D\bfv\w\ldots\w D^r\bfv)=\bfv\w D\bfv\w \ldots\w D^{r-1}\bfv\w D^{r+1}\bfv\\
&=&\bfv\w D\bfv\w \ldots\w (a_1(t)D^r\bfv-a_2(t)D^{r-1}\bfv+\ldots+(-1)^r\bfv)=\\
&=&a_1(t)\bfv\w D\bfv\w\ldots\w D^r\bfv=a_1(t)W_0(\bfv).
\end{eqnarray*}
The Wro\'nskian then takes the form (Abel's formula): \be W_0(\bfv)=W_0(\bfv)(0)\cdot \exp(\int_0^t
a(u)du),\label{eq:abelf} \ee where $W_0(\bfv)(0)$ denotes the value of the Wro\'nskian at $t=0$.
Equation~(\ref{eq:abelf}) shows that if $W(\bfv)(0)\neq 0$ then $W(\bfv)(t)\neq 0$ for all $t\in U$. We
shall see in Section~\ref{finalsec} why the proof of Liouville's theorem is a first example of the {\em
Schubert Calculus} formalism governing the intersection theory on Grassmann Schemes.

\claim{\bf Generalized Wronkians of Solutions of ODEs.}  Using {\em generalized Wro\'nskians} as
in~\ref{genwr}, Liouville's theorem~(\ref{eq:liouvthm}) can be rephrased as
\[
W_{(1)}(\bfv)=a_1(t)W_0(\bfv)\,,
\]
and generalized as follows.
\bclm{\bf Proposition.} {\em Let $1^k:=(1,1,\ldots,1)$ be the
{\em primitive partition} of the integer $1\leq k\leq r+1$. If $\bfv:=(v_0,v_1,\ldots, v_r)$ is a basis
of $\ker P_\bfa(D)$ then \be W_{(1^k)}(\bfv)=a_k(t)W(\bfv).\label{eq:genresliouv} \ee

} \eclm

Indeed, consider (\ref{eq:DifOp}). If $v\in \ker P_\bfa(D)$, then it is a $\KK$-linear combination of
$v_0,v_1,\ldots, v_r$, and hence the Wro\'nskian of these $r+2$ functions vanishes:
\[
W(v, v_0,v_1,\ldots, v_r)=0.
\]
By expanding the Wro\'nskian along the first column one obtains
\be
W(\bfv)D^{r+1}v-W_{(1)}(\bfv)D^rv+\ldots+(-1)^{r+1}W_{(1^{r+1})}(\bfv)v=0, \ee and combining with
$P_\bfa(D)=0$ this implies
\be \sum_{k=1}^{r+1}(-1)^k(W_{(1^k)}(\bfv)-a_k(t)W(\bfv))D^kv=0.\label{eq:wroncoef} \ee For general
$v\in\ker P_\bfa(D)$, the $(r+1)$-tuple $(v,Dv,\ldots, D^rv)$ is linearly independent, and
then~(\ref{eq:wroncoef}) implies~(\ref{eq:genresliouv}) for all $1\leq k\leq r+1$.

\claim{} \label{rmkpropor} A  natural question arises: Can we conclude that any generalized Wro\'nskian
$W_\blamb(\bfv)$ associated to a basis of $\ker P_\bfa(D)$ is a multiple of the Wro\'nskian $W_0(\bfv)$?
The answer is obviously yes. In fact  whenever one encounters one exterior factor in the generalized
Wro\'nskian  of the form $D^{j+\lambda_{r-j}}\bfv$ with $j+\lambda_{r-j}\geq r+1$, one  uses the
differential equation to express $D^{j+\lambda_{r-j}}\bfv$ as a linear combination of lower derivatives
of the vector $\bfv$, with  coefficients polynomial expressions in $\bfa$ and its derivatives,
\[
W_\blamb(\bfv)=G_\blamb(\bfa,D\bfa,D^2\bfa,\ldots)W(\bfv).
\]
The coefficient  $G_\blamb(\bfa,D\bfa,D^2\bfa,\ldots)$ assumes a particular interesting form in the case
the coefficients $\bfa$ of the equation are constant (so $D^i\bfa=0$, for $i>0$). We will address this
case in Section~\ref{finalsec}.

\section{Wro\'nski Sections of Line Bundles}\label{wrseclbclass}

{ \claim{}\label{bundles}  A {\em holomorphic vector bundle} of rank $d+1$ on a smooth complex
projective variety $X$  is a holomorphic map $\varrho:F\sra X$, where the complex manifold $F$ is
locally a product of $X$ and a complex $(d+1)$-dimensional vector space, cf.~~\cite[page 69]{GH}. For
$P\in X$, we denote by $F_P:=\varrho^{-1}(F)\subset F$ the fiber.

Consider the vector space $H^0(F):=H^0(X,F)$ of global holomorphic sections of $F$ (omitting the base
variety when clear from the context). For $s\in H^0(X,F)$ we will denote the value of $s$ at $P\in X$ by
$s(P)\in F_P$. The image of $s$ in the stalk of the {\em sheaf of sections of $F$} at $P$ will be
denoted by $s_P$.

 A {\em line bundle} over $X$ is a vector bundle of rank $1$. The set of isomorphism classes of line bundles
 on $X$ is a group under the tensor product; this group is denoted by $Pic(X)$.
 If $\pi:X\sra S$ is a proper flat morphism,  then we define a {\em relative line bundle}  as an equivalence
 class of line bundles on $X$,  where $\Lcal_1$ and
 $\Lcal_2$  are declared equivalent if $\Lcal_1\otimes \Lcal_2^{-1}\cong \pi^*\Ncal$, for some
 $\Ncal\in Pic(S)$. The group of isomorphism classes of  relative line bundles on
 $X$ is denoted by $Pic(X/S):=Pic(X)/\pi^*Pic(S)$.

\claim{}\label{charts} In the attempt to keep the paper self-contained, we recall a few basic notions
about line bundles  on a smooth projective complex curve.  From now on, we denote the curve by $C$. It
will often be identified with a {\em compact Riemann surface}, i.e. with  a complex manifold of complex
dimension $1$ equipped with a holomorphic atlas ${\frak A}:=\{(U_\alpha,z_\alpha)\,|\,
\alpha\in\Acal\}$, where $z_\alpha$ is a local coordinate on an open $U_\alpha$.  In this context,
denote by $\Ocal_C$ the sheaf of holomorphic functions on $C$: for $(U_\alpha,z_\alpha)\in {\frak A}$
the sheaf $\Ocal_C(U_\alpha)$ is the $\CC$-algebra of complex holomorphic functions in $z_\alpha$.

The {\em canonical line} bundle of $C$ is the line bundle $K\sra C$ whose transition functions
 are the derivatives of the coordinate changes,
\[\kappa_{\alpha\beta}:U_\alpha\cap U_\beta\lra \CC^*\,,  \ \ \kappa_{\alpha\beta}={dz_\alpha/
dz_\beta}\,.
\]
The holomorphic functions $\{\kappa_\ab\}$ obviously form a {\em cocycle}: $\kappa_\ab
\kappa_{\beta\gamma}=\kappa_{\alpha\gamma}$. A global holomorphic section $\omega\in H^0(C,K)$ is a {\em
global holomorphic differential}, i.e. a collection $\{f_\alpha\,dz_\alpha\}$, where $f_\alpha\in
\Ocal(U_\alpha)$ and ${f_\alpha}_{|_{U_\alpha\cap
U_\beta}}=\kappa_{\alpha\beta}{f_\beta}_{|_{U_\alpha\cap U_\beta}}$. We shall write
$f_\alpha\,dz_\alpha= \omega_{|_{U_\alpha}}$. The integer   $g=h^0(K):=\dim_\CC H^0(K)$ is
  the {\em genus of the curve}.

\claim{\bf Jets of line bundles.} \label{jetslb1} Let $\gf$ be a proper flat family of smooth projective
curves of genus $g\geq 1$ parameterized by some smooth scheme $S$. Let $\TS\times_S\TS\sra S$ be the
$2$-fold fiber product of $\TS$  over $S$ and let  $p,q:\TS\times_S \TS\sra \TS$ be the projections onto
the first and the second factor respectively. Denote by  $\delta:\TS\sra \TS\times_S \TS$ be the
diagonal morphism and by $\Ical$ the ideal sheaf of the diagonal in $\TS\times_S \TS$. The {\em relative
canonical bundle}  of  the family $\pi$ is by definition $K_\pi:=\delta^*(\Ical/\Ical^2)$.  For each
$\Lcal \in Pic(\TS/S)$, see \ref{bundles}, and  each $h\geq 0$ let
\be
J^h\Lcal:=p_*\left({O_{\TS\times _S\TS}\over \Ical^{h+1}}\otimes q^*\Lcal\right)\label{eq:defprpar}
\ee
 be the  bundle of {\em jets}  (or {\em principal parts}) of $\Lcal$ of order $h$. As $\TS$ is smooth,
 $J^h\Lcal$ is  a vector bundle on $\TS$ of rank $h+1$.

  By definition,  $J^0\Lcal=\Lcal$. Set, by convention, $J^{-1}\Lcal=0$ -- the vector bundle of rank $0$.
  The fiber  of $J^h\Lcal$ over $P\in \TS$ -- a complex vector space of dimension $h+1$ -- will be denoted
  by $J^h_P\Lcal$. The obvious exact sequence
\begin{center}
${\displaystyle 0\lra {\Ical^{h}\over\Ical^{h+1}}\lra {O_{\TS\times_S \TS}\over \Ical^{h+1}} \lra
{O_{\TS\times_S\TS}\over \Ical^{h}}\lra 0,}$
\end{center}
gives rise to an exact  sequence  (see~\cite[p.~224]{laksov} for details):
\be
0\lra \Lcal\otimes K_\pi^{h+1}{-\hspace{-4.5pt}\lra}
J^h\Lcal\stackrel{t_{h,h-1}}{-\hspace{-1.6pt}-\hspace{-4.6pt}\lra} J^{h-1}\Lcal\lra 0.\label{eq:1stexseq}
\ee
 If $\pi_0\,:\,C\sra \{pt\}$ is a trivial family over a point,  i.e. reduced to  a single curve,
 and if $L$ is any line bundle, then the exact sequence~(\ref{eq:1stexseq}) for $J^hL$ remains the same:
 in this case the relative canonical bundle  coincides with the canonical bundle of the curve.

\claim{} In notation of Section~\ref{charts}, let $v=(v_\alpha)$  be a non-zero holomorphic section of a
line bundle $L$, i.e. $v_\alpha\in \Ocal(U_\alpha)$ and $v_\alpha=\ell_\ab\cdot v_\beta$ on $\Uab$,
where $\{\ell_\ab\}$ are transition functions. Let $(U_\alpha,z_\alpha)$ be a coordinate chart of $C$
trivializing $L$. Denote by $D_\alpha:\Ocal(U_\alpha)\sra \Ocal(U_\alpha)$ the derivation $d/dz_\alpha$
and by $D_\alpha^j$ the $j$-th iterated of $D_\alpha$. Then
 \be
D_hv=\left\{\pmatrix{v_\alpha\cr D_\alpha v_\alpha\cr \vdots\cr D_\alpha^h v_\alpha} \,|\, \alpha\in
\Acal\right\}\label{eq:dhvsdhr}
\ee
is a section of $J^hL$ -- see~\cite{CuGaNi}. It may thought of as a {\em global derivative} of order $h$
of the section $v$. In fact it is a local representation of $v$ together with its first $h$ derivatives.

The truncation morphism occurring in~(\ref{eq:1stexseq}),  $t_{h,h-1}:J^hL\sra J^{h-1}L$, is defined in
such a way that $t_{h,h-1}(D_hv(P))=(D_{h-1}v)(P)$. See~\cite{CuGaNi} for further details.

\claim{}  One says that  $v\in H^0(L)$  vanishes at $P\in C$ with multiplicity at least $h+1$ if
$(\Dcal_hv)(P)=0$. Concretely, if $v_\alpha\in \Ocal_C(U_\alpha)$ is the local representation of $v$ in the open
set $U_\alpha$,  then $v$ vanishes at $P\in U_\alpha$ with multiplicity at least $h+1$ if $v_\alpha$ vanishes at
$P$ together with all of its first $h$ derivatives. The fact that $D_hv$ is a section of $J^hL$ says  that the definition of vanishing at a point $P$ does not depend on the open set $U_\alpha$ containing it.

We also say that the order of $v$ at $P$ is $h\geq 0$ if $D_{h-1}v(P)=0$ and $D_hv(P)\neq 0$. To each
$0\neq v\in H^0(L)$ one may attach a divisor on $C$: \be (v)=\sum_{P\in
C}(ord_Pv)P\,.\label{eq:finitesum} \ee The  sum~(\ref{eq:finitesum}) is finite because $v$ is locally a
holomorphic function and hence its zeros are isolated and the compactness of $C$ implies that they are
finitely many. The {\em degree} of $v$ is $\sum_{P\in C}ord_Pv\geq 0$. This number does not depend on a
holomorphic section of $L$, and  by definition is the {\em degree} of $L$. The degree of the canonical
bundle is $2g-2$ \cite[p.~8]{ACGH}. The set of isomorphism classes of line bundles of degree $d$ is
denoted by $Pic^d(C)$. If $\gf$ is a smooth proper family of smooth curves of genus $g$, then
$Pic^d(\TS/S)$ denotes the relative line bundles of relative degree $d$. A bundle $\Lcal\in Pic(\TS/S)$
has relative degree $d$ if $\deg(\Lcal_{|_{\TS_s}})=d$  for each $s\in S$.

\claim{} If $U$ is a (finite dimensional complex) vector space, $G(k,U)$ will denote the  Grassmannian
 parameterizing the $k$-dimensional vector subspaces of $U$.
 Let $g^r_d(L)$  be a point on chart $C$ of $G(r+1, H^0(L))$, where $L\in Pic^d(C)$.
 We write $g^r_d$  for $g^r_d(L)$ and some $L\in Pic^d(C)$.
 If $E=\sum e_PP$  is an effective divisor on $C$,  and $V$ is a $g^r_d(L)$, let
 \[
 V(-E):=\{v\in V\,|\, ord_Pv\geq e_P\},
 \]
Clearly $V(-E)$ is a vector subspace in $V$; it is not empty because it  contains at least the zero
section. If $\dim V(-P)=r$ for all $P\in C$, then the $g^r_d(L)$ is said to be {\em base point free}. It
is very ample if $\dim_\CC V(-P-Q)=r-1$ for all $(P,Q)\in C\times C$. If $V$ is base point free and
$\bfv:=(v_0,v_1,\ldots, v_r)$ is a basis of $V$, the map \be \left\{\matrix{\phi_\bfv&:&C&\lra&
\PP^r\cr\cr {}&{}&P&\longmapsto&(v_0(P): v_1(P):\ldots:v_r(P))} \right.\label{eq:Vmorph} \ee is a
morphism whose image is a projective algebraic curve of degree $d$. Although the complex value of  a
section at a point is not well defined, the ratio of two sections is. Thus the map~(\ref{eq:Vmorph}) is
well defined.  If $V$ is very ample, (\ref{eq:Vmorph}) is an embedding, i.e. a biholomorphism onto its
image.

\claim{}\label{canmorph} Let ${\bm\omega}:=(\omega_0,\omega_1,\ldots, \omega_{g-1})$ be a basis of $H^0(K)$. The
map
\[
\phi_{\bm\omega}:=(\omega_0:\omega_1:\ldots:\omega_{g-1}):C\sra \PP^{g-1}
\]
sending $P\mapsto (\omega_0(P):\omega_1(P):\ldots:\omega_{g-1}(P))$ is  the {\em canonical morphism},
that is, its image in $\PP^{g-1}$  is a curve of degree $2g-2$. If the  canonical morphism is not an
embedding, the curve is called {\em hyperelliptic}.

\bclm{\bf Definition.} {\em Let $V$ be a $g^r_d(L)$. A point  $P\in C$ is a $V$-ramification point if there exists $0\neq v\in V$ such that $D_rv(P)=0$, i.e. iff there exists a non-zero $v\in V$ vanishing at $P$ with multiplicity  $r+1$ at  least.}
\eclm
Ramification points of a $g^r_d$ can be detected as zero loci of suitable Wro\'nskians. Let
\[
\bfv:=(v_0, v_1,\ldots, v_r)
\]
be a basis of $V$ and let $v_{i,\alpha}:U_\alpha\sra \CC$ be holomorphic functions representing  the
restriction of the section $v_i$ to $U_\alpha$,  for $0\leq i\leq r$. If  $P\in U_\alpha$ is a
$V$-ramification point,  let $v=\sum_{i=0}^ra_iv_i$ be such that $D_rv(P)=0$. The last condition
translates into the following linear system: \be WM_\alpha(\bfv)\pmatrix{a_0\cr a_1\cr \vdots\cr
a_{r}}:=\pmatrix{v_{0,\alpha}&v_{1,\alpha}&\ldots&v_{r,\alpha}\cr D_\alpha v_{0,\alpha}&D_\alpha
v_{1,\alpha}&\ldots&D_\alpha v_{r,\alpha}\cr\vdots&\vdots&\ddots&\vdots\cr D_\alpha^r
v_{0,\alpha}&D_\alpha^r v_{1,\alpha}&\ldots&D_\alpha^r v_{r,\alpha}}\pmatrix{a_0\cr a_1\cr\vdots\cr
a_{r}}=\pmatrix{0\cr 0\cr\vdots\cr 0}.\label{eq:wmalfa} \ee It admits a non-trivial solution if and only
if the determinant
\[
W_0(\bfv_\alpha)=\bfv_\alpha\w D_\alpha\bfv_\alpha\w\ldots\w D_\alpha^{r}\bfv_\alpha\in \Ocal_C(U_\alpha)
\]
vanishes at $P$. It is easy to check that on $U_\alpha\cap U_\beta$ one has (see
e.g.~\cite[Ch.~2-18]{forster} or  \cite{CuGaNi})
\[
W_0(\bfv_\alpha)={\ell_{\alpha\beta}}^{r+1}(\kappa_{\alpha\beta})^{r(r-1)\over 2}W_0(\bfv_\beta),
\]
and thus the data $\{W_0(\bfv_\alpha)\,|\, \alpha\in \Acal\}$ glue  together to give a global holomorphic
section
\be
 W_0(\bfv)\in H^0(C, \wlbund),\label{eq:wrbund}
\ee said to be the {\em Wro\'nskian of  the basis} $\bfv$ of $V.$ The Wro\'nskian of any such a basis cannot vanish identically. Indeed, write the section $W_0(\bfv)$ as
\[
W_0(\bfv):= D_rv_0\w\ldots\w D_rv_r,
\]
where $D_rv$ is as in~(\ref{eq:dhvsdhr}), i.e. $D_rv_j$ is locally represented by the $j$-th row of the
matrix~(\ref{eq:wmalfa}). Assume that $W_0(\bfv_\alpha)$ vanished everywhere along the smooth connected
curve $C$. Then the sections $D_rv_j$, for $0\leq j\leq r$, corresponding to the columns of the
matrix~(\ref{eq:wmalfa}),  are linearly dependent, that is, up to a basis renumbering,
\[
D_rv_0=a_1D_rv_1+\ldots+a_rD_rv_r\in H^0(J^rL).
\]
However $D_r:H^0(L)\sra H^0(J^rL)$ is a section associated to the surjection $H^0(J^rL)\sra H^0(L)$,
induced by the truncation map $J^rL\sra L\sra 0$ (see e.g.~\cite[Section 2.7]{CuGaNi}) and one would
get the non-trivial linear relation
\[
v_0=a_1v_1+\ldots+a_rv_r\in H^0(L),
\]
against the assumption that $(v_0,v_1,\ldots,v_r)$ is  a basis of $V$.

As a consequence  the  ramification locus of the given $g^r_d$ occurs in codimension $1$. The construction
does not depend on the choice of a  basis $\bfv$ of $V$. Indeed,  if $\bfu$ were another one, then
$\bfu=A\bfv$ for some invertible $A\in Gl_{r+1}(\CC)$, and thence $ W_0(\bfu)=\det(A) W_0(\bfv)$. Thus
any basis of $V$ defines the same point of $\PP H^0(\wlbund)$, which we denote  by $W_0(V)$.

\claim\label{S39}{\bf The Wro\'nski Map.} We have so constructed   a  map: \be \left\{\matrix{G(r+1,
H^0(L))&\lra&\PP H^0(\wlbund)\cr\cr V&\longmapsto&W_0(V)}\right.\label{eq:wronskimap} \ee which
associates to each $g^r_d(L)$ its Wro\'nskian $W_0(V)$. Adopting the  same terminology used in the
literature when $C=\PP^1$ and $L:=O_{\PP^1}(d)$ (see e.g.~\cite{EreGab1},~\cite{EreGab2}),  the
map~(\ref{eq:wronskimap}) will be called  {\em Wro\'nski map}. Its behavior depends  on the curve and on
the  choice of the linear system.
 It is, in general, neither injective nor surjective as the following  two extremal  cases show. If $C=\PP^1$,
 the unique bundle of degree $d$ is $O_{\PP^1}(d)$, $K=O_{\PP^1}(-2)$ and the
 the Wro\'nski map
 \[
 G(r+1, H^0(O_{\PP^1}(d)))\lra \PP H^0(O_{\PP^1}((r+1)(d-r))),
 \]
in this case  defined between two varieties of the same dimension,  is a finite surjective morphism  of
degree equal to the Pl\"ucker degree of the Grassmannian $G(r+1,d+1)$. In particular it is not injective
-- see~\cite{EHcuspidal, Scherb2} and \cite{EreGab1,EreGab2}  over the real numbers.
  At a general point of $ \PP H^0(\Ocal_{\PP^1}((r+1)(d-r)))$ (represented by a form of degree $(r+1)(d-r)$) there correspond as many distinct  linear systems $V$ as the degree of the Grassmannian. For a closer analysis of the fibers of such a morphism see~\cite{Scherb1}.

  On the other hand if   $C$ is  hyperelliptic and $\Mcal\in Pic^2(C)$ is the line bundle defining its unique $g^1_2$, cf. Section~\ref{canmorph},  then $G(2, H^0(\Mcal))$ is just a point and the Wro\'nski map:
  \[
 G(2, H^0(\Mcal))\sra  \PP H^0(\Mcal^{\otimes 2}\otimes K)
  \]
is trivially injective and not surjective, as by Riemann-Roch formula $h^0(\Mcal^{\otimes 2}\otimes K)
>1$.

Later on we shall see how to make the situation more uniform, by enlarging in a natural way the notion
of linear system on a curve.  It will be  one of the bridges connecting this part of the survey with the
first one, regarding Wro\'nskians of differential equations.

\claim{\bf The $V$-weight of a point.}  Let $V$ be a $g^r_d$ and  $P\in C$. The {\em $V$-weight}
at $P$ is:
 \[
 wt_V(P):=ord_PW_0(V)=ord_PW_0(\bfv),
 \]
 for some basis $\bfv$ of $V$. The {\em total weight} of the $V$-ramification points is:
\[
wt_V=\sum_{P\in C}wt_V(P),
\]
where the above sum is clearly finite.
The {\em total weight} coincides with the degree of the bundle $\wlbund$, i.e. the degree of its first Chern
class:
\begin{eqnarray}
\wt_V&=&\int_Cc_1(\wlbund)\cap[C]=\nonumber\\
&=&(r+1)\int_C(c_1(L)\cap[C])+{r(r+1)\over 2}\int c_1(K)\cap [C]=\nonumber\\
&=&(r+1)d+(g-1)r(r+1),\label{eq:brillsegre}
\end{eqnarray}
which is  the so-called {\em Brill--Segre formula}. For example,  a smooth plane curve of degree $d$ can be
thought of as an abstract curve (compact Riemann surface) embedded in $\PP^2$ via some $V\in G(3, H^0(L))$ for
some $L\in Pic^d(C)$:
\[
(v_0:v_1:v_2):C\lra \PP^2
\]
where $\bfv:=(v_0,v_1,v_2)$ is a basis of $V$. The $V$-ramification points correspond, in this case,  to
flexes of the image of $C$ in $\PP^2$. According to the genus-degree formulae, the total number of
flexes, keeping multiplicities into account, is given by~(\ref{eq:brillsegre})  for $r=2$
\[
f=3d(d-2),
\]
which is one of the famous {\em Pl\"ucker formulas} for plane curves. \claim{\bf Wro\'nskians on
Gorenstein Curves.} Let $C$ be an irreducible  plane curve of degree $d$ with $\delta$ nodes and
$\kappa$ cusps. Using the extension of the Wro\'nskian of a linear system defined on a Gorenstein curve,
due to Widland and Lax \cite{WL1}, the celebrated Pl\"ucker formula
\[
f=3d(d-2)-6\delta-8\kappa
\]
can be obtained from the tautological identity (see~\cite{Gat0} for details):
\begin{center}
$\sharp$(smooth $V$--ramification points) =\\
=$\sharp$(ramification points) - $\sharp$(singular ramification points).
\end{center}
For more on jets and Wro\'nskians on Gorenstein curves see~\cite{EST1} and~\cite{EST2}.

\claim{}\label{orderpartition} The $V$-weight of a point $P$ coincides with the weight of its {\em order
partition}. We say that $n\in\NN$ is a $V$-order at $P\in C$ if there exists $v\in V$ such that
$ord_Pv=n$. Each point possesses only $r+1$ distinct $V$-orders. In  fact $n$ is a $V$-order if $\dim
V(-nP)>\dim V(-(n+1)P)$. We have the following sequence of inequalities:
\begin{center}
$r+1=\dim V\geq \dim V(-P)\geq \dim V(-2P)\geq \ldots\geq \dim V(-dP)\geq \dim V(-(d+1)P)=0$
\end{center}
The last dimension is zero because the unique section of $V$ vanishing at $P$ with multiplicity $d+1$ is
zero. At each step the dimension does not drop more than one unit and then there must be precisely $r+1$
jumps. If
\[
0\leq i_0<i_1<\ldots <i_r\leq d
\]
is the order sequence at some $P\in C$,  the $V$-{\em order partition at $P$} is
\[
\blamb(P)=(i_r-r,i_{r-1}-(r-1),\ldots, i_1-1,i_0).
\]
One may choose a basis $(v_0,v_1,\ldots, v_r)$ of $V$ such that $ord_Pv_j=i_j$. The use of such  a basis
shows that the Wro\'nskian $W_0(\bfv)$ vanishes at $P$ with multiplicity
\begin{center}
$\displaystyle{wt_V(P)=\sum_{j=0}^r(i_j-j)=|\blamb(P)|}$
\end{center}
The following result is due to~\cite{Ponza} (unpublished) and to ~\cite{towse}.
\bclm{\bf Proposition.} \label{potow} {\em Partition $\blamb$ is the $V$-order partition at
$P\in C$ if and only if $W_\bmu(\bfv_\alpha)(P)=0$, for all $\bmu\neq \blamb$ such that $|\bmu|\leq
|\blamb|$, and $W_\blamb(\bfv_\alpha)(P)\neq 0$ (here $\bfv_\alpha$ is any local representation of a
basis of $V$ around $P$). }
\eclm
In this case  $W_0(V)$ vanishes at $P$ with multiplicity exactly $|\blamb|$.

\claim{} \label{intrwrsec} A more intrinsic way to look at Wro\'nskians and ramification points, which
can be  generalized to the case of families of curves, is as follows. For $V\in G(r+1, H^0(L))$ one
considers the vector bundle map \be \Dcal_r:C\times V\lra J^rL\label{eq:vbmsr} \ee defined by
$\Dcal_r(P,v)=D_rv(P)\in J^r_PL$. Both bundles have rank $r+1$ and since $V$ has only finitely many
ramification points, there is a non-empty open subset of $C$ where the map $\Dcal_r$ has the maximal
rank $r+1.$ Then $P\in C$ is a $V$-ramification point if  $\rk_P\Dcal_r\leq r.$ The rank of $\Dcal_r$ is
smaller than the maximum if and only if the determinant map of~(\ref{eq:vbmsr})
\begin{center}
$\bw^{r+1}\Dcal_r:O_C\lra \bw^{r+1}J^rL$
\end{center}
vanishes at $P.$ The section $\bw^{r+1}\Dcal_r\in H^0(\bw^{r+1}J^rL)=H^0(\wlbund)$ is precisely the {\em
Wro\'nski section}, which vanishes precisely where the map $\Dcal_r$ drops rank. If
$v=a_0v_0+\ldots+a_rv_r$, with respect to the basis $\bfv=(v_0,v_1,\ldots, v_r)$ of $V$,  then
$D_rv=a_0D_rv_0+a_1D_rv_1+\ldots+a_rD_rv_r.$ On a trivializing open set $U_\alpha$ of $C$ one has the
expression:
\[
(D_rv)_{|_{U_\alpha}}=\pmatrix{a_0v_{0,\alpha}+a_1v_{1,\alpha}+\ldots+a_rv_{r,\alpha}\cr a_0D_\alpha
v_{0,\alpha}+a_1D_\alpha v_{1,\alpha}+\ldots+a_rD_\alpha v_{r,\alpha}\cr \vdots\cr a_0D^r_\alpha
v_{0,\alpha}+a_1D^r_\alpha v_{1,\alpha}+\ldots+a_rD^r_\alpha v_{r,\alpha}}=W_0(\bfv_\alpha)\cdot
\pmatrix{a_0\cr a_1\cr\vdots\cr a_r}.
\]
In other words, the local representation of the map $\Dcal_r$ is:
\[
W_{0}(\bfv_\alpha):U_\alpha\times \CC^{r+1}\lra U_\alpha\times\CC^{r+1}
\]
from which:
\[
\det({\Dcal_r}_{|_{U_\alpha}})=\bfv_\alpha\w D_\alpha\bfv_\alpha\w\ldots\w D^r_\alpha\bfv_\alpha,
\]
i.e. $\bw^{r+1}\Dcal_r$ is represented by  the {\em Wro\'nskian} $W_0(\bfv).$ Changing the basis $\bfv$
of $V$, the Wro\'nski section gets multiplied by a non-zero complex number and hence:
\[
\bw^{r+1}\Dcal_r\,\,\mod\, \CC^*=W_0(V)\in \PP H^0(\wlbund)
\]
i.e. precisely the Wro\'nskian associated to the linear system $V.$

\claim{} \label{relativb} How do generalized Wro\'nskians come into play in this picture? Here the
question is more delicate. We have already mentioned that if the $V$-order partition of a point $P$ is
$\blamb(P)$ then the generalized Wro\'nskians $W_\bmu(V)$ must vanish for all $\bmu$ such that
$|\bmu|<|\blamb(P)|$ and $W_\blamb(P)\neq 0$. It is however well known that the general $g^r_d$ on a
general curve $C$ has only simple ramification points, i.e. all the points have weight $1$. This says
that if a $g^r_d$ has a ramification point with weight bigger than $1$, the generalized Wro\'nskians do
not impose independent conditions, as the locus occurs in codimension $1$ while the expected codimension
is bigger than $1$.

To look for more  geometrical content one can  move along two
directions. The first, that we just sketch here, consists in considering families of curves.

 Let $\gf$ be a proper flat family of smooth curves of genus $g$ and let $(\Vcal, \Lcal)$ be a relative $g^r_d$, i.e. $\Vcal$ is a locally free subsheaf of $\pi_*\Lcal$ and $\Lcal\in Pic^d(\TS/S)$.
One can  then study the ramification locus of the relative $g^r_d$ which fiberwise cuts the ramification locus of
$\Vcal_s\in G(r+1, H^0(\Lcal_{|_{\TS_s}}))$ through the degeneracy locus of the map
\[
\Dcal_r:\pi^*\Vcal\lra J^r_\pi\Lcal,
\]
where $J^r_\pi\Lcal$ denotes the jets of $\Lcal$ {\em along the fibers} (see e.g.~\cite{GatPon}). The
map above induces a section $\Ocal_\TS\sra \bw^{r+1}J^r_\pi\Lcal\otimes\bw^{r+1}\pi^*\Vcal$, which is
the relative Wro\'nskian $W_0(\Vcal)$ of the family.  Because of the exact sequence~(\ref{eq:1stexseq}):
\[
\bw^{r+1}J^r_\pi\Lcal\otimes\bw^{r+1}\pi^*\Vcal=\Lcal^{\otimes r+1}\otimes K_\pi^{\otimes {r(r+1)\over
2}}\otimes\bw^{r+1}\pi^*\Vcal\,.
\]
In this case the class in $A_*(\TS)$ of the ramification locus of $\Vcal$ is
\[
[Z(W_0(\Vcal))]=c_1(\Lcal^{\otimes r+1}\otimes K_\pi^{\otimes {r(r+1)\over 2}})-\pi^*c_1(\Vcal)\,.
\]

A  second approach to enrich the phenomenology of ramification points  consists in keeping the curve
fixed and varying the linear system. This is the only possible approach with curves of genus $0$: all
the smooth rational curves are isomorphic, and all the $g^r_d$s, with base points or not,  are
parameterized by the Grassmannian $G(r+1, H^0(O_{\PP1}ÿ(d)))$. Here the situation is as nice as one
would desire:  all what may potentially occur it occurs indeed. For instance,  if $\blamb_1,\ldots,
\blamb_h$ are partitions such that $\sum |\blamb_i|=(r+1)(d-r)$ (= the total weight of the ramification
points of a $g^r_d$) and $P_1,\ldots, P_h$ are arbitrary points on $\PP^1$ one can count the number of
all of the linear system such that the $V$ order partition at $P_i$ is precisely $\blamb_i$. However if
$C$ has higher genus, such a  kind of analysis is not possible anymore. For instance the general curve
$C$ of genus $g\geq 2$ has only simple Weierstrass points, i.e. all have weight $1$, but each curve
carries one and only one canonical system. The picture holding for linear systems on the  projective
line can be generalized in the case of higher genus curves provided one updates the notion of $g^r_d(L)$
to that of {\em a section of a Grassmann bundle}, a  path  which was first indicated
in~\cite{GatSalehyan} and then further developed in~\cite{NigroTesi} and~\cite{CuGaNi}.  Go to the next
two sections for a sketch of the construction.

\section{Wro\'nskians of Sections of Grassmann Bundles (in general)}\label{wrsecbngen}
This section is a survey of the contruction appeared in~\cite{NigroTesi}, partly published in~\cite{CuGaNi},
 with some applications
in~\cite{GatSalehyan}.

\claim{}
 Let $\varrho_d:F\sra X$ be a vector bundle of rank $d+1$ over a smooth complex projective variety $X$ of dimension $m\geq 0$. For each $0\leq r\leq d$, let $\varrho_{r,d}: G(r+1, F)\sra X$ be the Grassmann bundle of $(r+1)$-dimensional subspaces of the fibers of $F$. For  $r=0$ we shall write $\varrho_{0,d}:\PP(F)\lra X$,  where  $\PP(F):=G(1, F)$ is the {\em projective bundle} associated to $F$.
 The bundle $G(r+1, F)$ carries  universal exact sequence (cf.~\cite[Appendix B.5.7]{Fu1}):
\be
0\lra\Scal_r \stackrel{\iota_r}{\lra} \varrho_{r,d}^*F\lra \Qcal_r\lra 0,\label{eq:univexseq}
\ee where $\Scal_r$ is the {\em universal subbundle} of $\varrho_{r,d}^*F$ and $\Qcal_r$ is the {\em universal quotient bundle}.

Let
\begin{center}
$\Gamma(\varrho_{r,d}):=\{$holomorphic $\gamma:X\sra G(r+1, F)\,|\, \varrho_{r,d}\circ\gamma=id_X\}$
\end{center}
be the set of holomorphic sections of $\varrho_{r,d}$. The choice of  $\gamma\in \Grd$ amounts  to
specify a vector sub-bundle of $F$ of rank $r+1$. In fact the pull-back $\gamma^*\Scal_r$ via
$\gamma\in \Grd$ is  a rank $r+1$  subbundle of $F$. Conversely, given  a rank $r+1$ subbundle $\Vcal$
of $F$, one may define the section $\gamma_\Vcal\in\Grd$ by $\gamma_\Vcal(P)=\Vcal_P\in G(r+1, F_P)$.
The set $\Grd$ is huge and may have a very nasty behavior: even the case when $X=\PP^1$ and
$F=J^dO_{\PP^1}(d)$,  is far from being trivial. In fact  it is related with the small quantum
cohomology of Grassmannians, see~\cite{Bertram}.  A first simplification is to fix  $\xi\in Pic(X)$ to
study the space
\begin{center}
$\Gamma_\xi(\varrho_{r,d})=\{\gamma\in\Grd\,|\, \bw^{r+1}\gamma^*\Scal_r=\xi\}$.
\end{center}
Again, if $\xi=O_{\PP^1}(n)$ and $F=J^dO_{\PP^1}(d)$, then $\Gamma_{n}(\varrho_{r,d}):=\Gamma_{O_{\PP^1}(n)}(\varrho_{r,d})$ can be identified with the space of the holomorphic maps $\PP^1\sra G(r+1, d+1)$ of degree $n$, compactified in~\cite{Bertram} via a Quot-scheme construction. We shall see the easiest case ($n=0$) in Section~\ref{wmapprojline}. In the following, for our limited purposes, we shall restrict the attention to the  definitely
simpler set
\begin{center}
$\Gtrdb:=\{\gamma\in\Grd\,|\, \gamma^*\Scal_r$ is a Êtrivial rank $(r+1)$ subbundle of $F\sra X\}$.
\end{center}

 \bclm{\bf Proposition.} \label{opensubset} {\em The set $\Gtrdb$, if non empty,
can be identified with an open set of the Grassmannian $G(r+1, H^0(F))$. } \eclm \proof If $\gamma\in
\Gtrdb$, there is an isomorphism  $\phi:X\times \CC^{r+1}\sra \gamma^*\Scal_r$. Then
$\psi:=\gamma^*(\iota_r)\circ \phi: X\times \CC^{r+1}\sra F$ is a bundle monomorphism. Let
$\sigma_i:X\sra F$ defined by $\sigma_i(P)=\psi(P, \bfe_i)$. It is clearly a holomorphic section of $F$.

Furthermore $\sigma_0,\sigma_1,\ldots,\sigma_r$ span an $(r+1)$-dimensional subspace $U_\gamma$ of
$H^0(F)$ which does not depend on the choice of the isomorphism $\phi$. Thus $\gamma^*\Scal_r$ is
isomorphic to $X\times U_\gamma$ and $\gamma(P)=\{u(P)\,|\, u\in U_\gamma\}\in G(r+1, F_P)$. Conversely,
if $U\in G(r+1, H^0(F))$, one constructs a vector bundle morphism $\phi: X\times U\sra F$ via $(P,
u)\mapsto u(P)$. This morphism drops rank if $\bw^{r+1}\phi=0$,  this is a closed condition and so there is
 an open set  $\Ucal\subseteq G(r+1, F)$ such that for $U\in \Ucal$, the map $\phi_U$ makes $X\times
U$ into a vector subbundle of $F$. One so obtains a section $\gamma_U$ by setting
$\gamma_U(P)=U_P\in G(r+1, F_P)$. The easy  check that $\gamma_{U_\gamma}=\gamma$ and that
$U_{\gamma_U}=U$ is left to the reader.\qed

\claim{} \label{filtration}
 Assume now that $F$ comes equipped with a system  $F_\bullet$  of bundle epimorphisms $q_{ij}:F_i\lra
F_j $,  for each   $-1\leq j\leq i\leq d$, such that $F_d=F$, where  $F_i$ has rank $i+1$,
$q_{ii}=id_{F_i}$ and $q_{ij}q_{jk}=q_{ik}$ for each triple $d\geq i\geq j\geq k\geq -1$. We  set
$F_{-1}=0$ by convention. The map $q_{dj}:F\sra F_j$ will be simply denoted by $q_j$ and
$\{\ker(q_{i})\}$ gives a filtration of $F$ by subbundles of rank $d-i$. Let
\[
\partial_i:\Scal_r\lra \varrho_{r,d}^*F_i
\]
be the composition of the universal monomorphism $\Scal_r\sra \varrho_{r,d}^*F$ with the map $q_i$. The
universal morphism $\iota_r$  can be so identified with $\partial_d$.

\claim{} \label{svf44} For each $\blamb\in \Pcal^{(r+1)\times(d-r)}$ the subscheme:
 \be
 \Omega_\blamb(\varrho_{r,d}^*F_\bullet)=\{\Lambda\in G(r+1, F)\,|\,
\rk_\Lambda\partial_{j+\lambda_{r-j}-1}\leq j, \qquad 0\leq j\leq r\},\label{eq:omlamb}
\ee
of $G(r+1, F)$, is the
$\blamb$-{\em Schubert variety} associated to the system  $F_\bullet$ and to the partition $\blamb$.
The Chow classes modulo rational equivalence $\{[\Omega_\blamb(\varrho_{r,d}^* F_\bullet)]\,|\,\blamb\in \Pcal^{(r+1)\times (d-r)}\}$ freely generate
$A_*(G(r+1, F))$ as a module over $A^*(X)$ through the structural map  $\varrho_{r,d}^*$.

\claim{}\label{defNh} For each
$0\leq h\leq d+1$, let $N_h(F):=\ker (F\stackrel{q_{d-h}}{-\hskip-6pt-\hskip-6pt\lra} F_{d-h})$. It is a vector bundle of
rank $h$. One can define Schubert varieties according to such a {\em kernel} flag $N_\bullet(F)$ by
setting, for each partition $\blamb$ of length at most $r+1$:
\[
\Omega_\blamb(\varrho_{r,d}^*N_\bullet(F))=\{\Lambda\in G(r+1, F)\,|\, \Lambda\cap
N_{d+1-(j+\lambda_{r-j})}(F)\geq r+1-j\}\,.
\]
It is a simple exercise of linear algebra to show that
\[
\Omega_\blamb(\varrho_{r,d}^*F_\bullet)=\Omega_\blamb(\varrho_{r,d}^*N_\bullet(F)).
\]
Both descriptions are useful according to the purposes. The first description is more suited to describe
Weierstrass points as in Section~\ref{wrseclbclass}  (it gives an  algebraic generalization of the rank
sequence in a Brill-N\"other matrix, see~\cite[p.~154]{ACGH}), while the second is useful when dealing
with linear systems on the projective line (see Section~\ref{wmapprojline} below).

\bclm{\bf Definition.}\label{def23} {\em The $F_\bullet$-{\em Wro\'nskian} subvariety of $G(r+1, F)$ is
\[
{\frak W}_0(\varrho_{r,d}^* F_\bullet):=\Omega_{(1)}(\varrho_{r,d}^* F_\bullet).
\]
}
\eclm

By~(\ref{eq:omlamb}), the $F_\bullet $-Wro\'nski variety ${\frak W}_0(\varrho_{r,d}^* F_\bullet)$ of
$G(r+1, F)$ is the degeneracy scheme of the natural map $\partial_r:\Scal_r\lra \varrho_{r,d}^*F_r$,
i.e. the zero scheme of the map
\[
\bw^{r+1}\partial_r:\bw^{r+1}\Scal_r\lra \bw^{r+1}\varrho_{r,d}^*F_r.
\]

 The map \be W_0(\varrho_{r,d}^*F_\bullet):=\bw^{r+1}\partial_r\in
Hom(\bw^{r+1}\varrho_{r,d}^*\Scal_r, \bw^{r+1}\varrho_{r,d}^*F_r)=H^0(X,
\bw^{r+1}\varrho_{r,d}^*F_r\otimes\bw^{r+1}\varrho_{r,d}^*\Scal_r^\vee),\label{eq:wrsecvb} \ee is the
{\em Wro\'nski section} (of the line bundle
$\bw^{r+1}\varrho_{r,d}^*F_r\otimes\bw^{r+1}\varrho_{r,d}^*\Scal_r^\vee$). The $F_\bullet$-Wro\'nski
variety is then a Cartier divisor, because  it is the zero scheme of the Wro\'nski
section~(\ref{eq:wrsecvb}). In this setting, the  Schubert subvariety
$\Omega_\blamb(\varrho_{r,d}^*F_\bullet)$ of  $G(r+1, F)$, associated to the partition $\blamb\in
\Pcal^{(r+1)\times(d-r)}$, plays the role of a generalized Wro\'nski subvariety associated to the system
$F_\bullet$.

\claim{} \label{bsptlc} Among all such Schubert varieties associated to $F_\bullet$  one can recognize
some distinguished ones. It is natural to define  the  $F_\bullet$-{\em base locus} subvariety of
$G(r+1, F)$  as
\[
\Bcal(\varrho_{r,d}^*F_\bullet)=\Omega_{(1^{r+1})}(\varrho_{r,d}^*F_\bullet);
\]
and  the  $F_\bullet$-{\em cuspidal locus} subvariety as
\[
\Ccal(\varrho_{r,d}^*F_\bullet)={\Omega}_{(1^r)}(\varrho_{r,d}^*F_\bullet).
\]
Each Schubert  subvariety ${\Omega}_\blamb(\varrho_{r,d}^*F_\bullet)$ has codimension $|\blamb|$ in
$G(r+1, F)$. In particular, the base locus variety $\Bcal(\varrho_{r,d}^*F_\bullet)$ has codimension
$r+1$.

\claim{} Let $\gamma\in \Grd$. The $F_\bullet$-{\em ramification locus} of $\gamma$  is the subscheme
$\gamma^{-1}({\frak W}_0(\varrho_{r,d}^*F_\bullet))$ of $X$,   its $F_\bullet$-{\em base locus}  is
$\gamma^{-1}(\Bcal(\varrho_{r,d}^*F_\bullet))$ and its $F_\bullet$-{\em cuspidal locus}  is
$\gamma^{-1}(\Ccal(\varrho_{r,d}^*F_\bullet))$. The definition of {\em Wro\'nski map} defined on
sections of Grassmann bundles equipped with filtrations,
 as in Section~\ref{filtration},  is very natural too.

\bclm{\bf Definition.} {\em For $\gamma\in \Gamma(\varrho_{r,d})$, the section

\begin{center}
$\displaystyle{W_0(\gamma):=\gamma^*(W_0(\varrho_{r,d}^*F_\bullet))\in
H^0(X,\bw^{r+1}F_r\otimes\bw^{r+1}\gamma^*\Scal_r^\vee)}$
\end{center}
will be called the $F_\bullet$-{\em  Wro\'nskian} of $\gamma$. } \eclm The class  in $A^*(X)$ of the
ramification locus of $\gamma$ is:

\begin{eqnarray}
[Z(W_0(\gamma))]&=&
[\gamma^{-1}({\frak W}_{0}(\varrho_{r,d}^*F_\bullet))]=\gamma^*[{\frak W}_0(\varrho_{r,d}^*F_\bullet)]
=\nonumber\\
&=&c_1(\bw^{r+1}F_r\otimes\bw^{r+1}\gamma^*\Scal_r^\vee)\cap[X]= (c_1(F_r)-
\gamma^*c_1(\Scal_r))\cap[X].\label{eq:classramloc}
\end{eqnarray}
If $X$ is a curve, the expected dimension of the ramification locus is $0$ and so, when $\gamma$ is not
entirely contained in the Wro\'nski variety,  the {\em total weight} $w_\gamma$ of the ramification
points of $\gamma$ is by definition the degree of the cycle $[\gamma^{-1}({\frak
W}_0(\varrho_{r,d}^*F_\bullet))]$:
\[
w_\gamma=\int_X(c_1(F_r)-\gamma^*c_1(\Scal))\cap[X].
\]
According to the definitions above, a point $P\in X$ is a ramification point of $\gamma\in
\Gamma(\varrho_{r,d})$ if $W_0(\gamma)(P)=0$,  which amounts to say  that the map
$\gamma^*\partial_r: \gamma^*\Scal_r\sra F_r$ drops rank at $P$.

\bclm{\bf Definition.} {\em Fix $\xi\in Pic(X)$. The holomorphic  map:
\[
\left\{\matrix{\Gamma_\xi(\varrho_{r,d})&\lra&\PP H^0(\bw^{r+1}F_r\otimes
\xi^\vee)\cr\cr\gamma&\longmapsto&W_0(\gamma)\,\,\pmod {\CC^*}}\right.
\]
 is the  {\em Wro\'nski map} defined on $\Gamma_\xi(\varrho_{r,d})$.
 }
\eclm

Indeed $W_0(\gamma)$ is a section of
$\gamma^*(\bw^{r+1}\varrho_{r,d}^*F_r\otimes\bw^{r+1}\Scal_r^\vee)=\bw^{r+1}F_r\otimes \xiÿ^\vee$.
The class of the {\em
ramification locus} of $\gamma$, as in~(\ref{eq:classramloc}), can be now expressed as:
\begin{center}
$[Z(W_0(\gamma))]= (c_1(F_r)-\xi))\cap [X]\in A^*(X)$.
\end{center}
\claim{\bf The  Extended Wro\'nski Map.} It is particularly easy  to express the Wro\'nskian of a
section $\gamma\in \Gamma_{\tt triv}(\varrho_{r,d})$. Let $U\in G(r+1, F)$ such that $\gamma=\gamma_U$.
The pull-back of the map $
\partial_{r}:\Scal_r\lra \varrho_{r,d}^*F_r
$ is
\be
\gamma^*\partial_{r}:X\times U\lra F_r.\label{eq:pllbck}
\ee
The Wro\'nskian is the determinant of the  map~(\ref{eq:pllbck}):
\[
\bw^{r+1}\gamma^*\partial_{r}:\bw^{r+1}(X\times U)\sra \bw^{r+1}F_r.
\]
Once a basis $(u_0,u_1,\ldots, u_r)$ of $U$ is chosen, the Wro\'nskian
\[
\bw^{r+1}\gamma^*\partial_r\in H^0(X, \bw^{r+1}F_r)
\]
is represented by the holomorphic section $X\sra \bw^{r+1} F_r$ given by:

\begin{center}
$P\mapsto q_r(u_0)(P)\w q_r(u_1)(P)\w\ldots \w q_r(u_r)(P)\in \bw^{r+1}F_P,$
\end{center}
where $q_r$  is the epimorphism introduced  in~\ref{filtration}.
 Changing basis the section gets multiplied by a non-zero constant, and so the Wro\'nski  map
\begin{center}
$\Gamma_{\tt triv}(\varrho_{r,d})\lra \PP H^0(X, \bw^{r+1}F_r)$
\end{center}
defined by $\gamma\mapsto W_0(\gamma)\,\mod\,\,\CC^* \in \PP H^0(X, \bw^{r+1}F_r)$ coincides with the map
\[
\left\{\matrix{G(r+1, H^0(F))&\lra&\PP H^0(X, \bw^{r+1}F_r)\cr\cr U&\longmapsto& q_r(u_0)\w q_r(u_1)\w\ldots\w
q_r(u_r)\,\, \mod \,\, \CC^*}\right.
\]
where $\bfu=(u_0,u_1,\ldots, u_r)$ is any basis of $U$.

 \claim{} \label{intthgrbund} Here is a quick review of intersection theory on $G(r+1, F)$
  which is necessary for enumerative geometry purposes. First recall some basic terminology and notation.
  Let $a=a(t)=\sum_{n\geq 0}a_nt^n$ be a formal power series with coefficients in some ring
  $A$ and $\blamb$ be a partition as in~(\ref{eq:lambpart}). Set  $a_n=0$ for $n<0$.
  The $\blamb$-{\em Schur polynomial} associated to $a$  is, by definition:
 \be
\Delta_\blamb(a)=\det(a_{i+\lambda_{r-i}-j})_{0\leq i,j\leq r}
=\left|\matrix{a_{\lambda_r}&a_{\lambda_{r-1}+1}&\ldots&a_{\lambda_0+r}\cr
a_{\lambda_r-1}&a_{\lambda_{r-1}}&\ldots&a_{\lambda_0+r-1}\cr \vdots&\vdots&\ddots&\vdots\cr
a_{\lambda_r-r}&a_{\lambda_{r-1}-(r-1)}&\ldots&a_{\lambda_0}}\right|\,.\label{eq:defschur}
 \ee

The Chern polynomial of a bundle $\Ecal$ is denoted by $c_t(\Ecal)$. Write
$c_t(\Qcal_r-\varrho_{r,d}^*F))$ for the ratio $c_t(\Qcal_r)/c_t(\varrho_{r,d}^*F)$ of Chern
polynomials. According to the  Basis Theorem \cite[p.~268]{Fu1}, the Chow group $A^*(G(r+1, F))$ is a free
$A^*(X)$-module (via the structural morphism $\varrho_{r,d}^*:A^*(X)\sra A^*(G(r+1, F))$) generated by
 \[
 \{\,\Delta_\blamb(c_t(\Qcal_r-\varrho_{r,d}^*F))\cap [G(r+1, F)]\,|\, \blamb\in \Pcal^{(r+1)\times (d-r)}\,\},
 \]

 If $r=0$, let
 \[
 \mu^i:=(-1)^ic_1(\Scal_0)^i\cap [\PP(F)]
 \]
 for each $i\geq 0$.
 Then, by~\cite[Ch.~14]{Fu1},
 $(\mu^0,\mu^1,\ldots, \mu^d)$ is an $A^*(X)$-basis of $A^*(\PP(F))$ and for each $j\geq 0$ the following relation, defining the Chern classes of $F$, holds:
 \be
  \mu^{d+1+j}+\varrho_{0,d}^*c_1(F)\mu^{d+j}+\ldots + \varrho_{0,d}^*c_{d+1}(F)\mu^j=0\,.\label{eq:relchowpf}
 \ee
A main result of  \cite{GatSant1} says that  $\bw^{r+1}A^*(\PP(F))$ can be equipped with a structure of
$A^*(G(r+1, F))$-module of rank $1$. It is generated by $\mu^0\w\mu^1\w\ldots\w\mu^r$ in such a way
that, for each $\blamb\in \Pcal^{(r+1)\times (d-r)}$,
 \be
\Delta_\blamb(c_t(\Qcal_r-\varrho_{r,d}^*F))\cdot \mu^0\w\mu^1\w\ldots\w\mu^r=
\partoformu.\label{eq:gatsant}
 \ee
We shall see in the last section that  $\Delta_\blamb(c_t(\Qcal_r-\varrho_{r,d}^*F))$ are related to
Wro\'nskians associated to a fundamental system of solutions of a suitable differential equation. Define
now:
\[
\ep^{i}:=[\Omega_{(i)}(\varrho_{0,d}^*F_\bullet)]\in A^*(\PP(F)),\quad 0\leq i\leq d,
\]
where $\Omega_{(i)}(\varrho_{0,d}^*F_\bullet)$ is nothing but the  zero locus in codimension $i$ of the map
$\partial_{i-1}:\Scal_0\sra F_{i-1}$.  Because of the relation:
\be
\ep^i=\sum_{j=0}^i\varrho_{0,d}^*c_j(F_{i-1})\mu^{i-j},\label{eq:relmuep}
\ee
 it follows that $(\ep^0,\ep^1,\ldots, \ep^d)$ is an $A^*(X)$-basis of $A^*(\PP(F))$ as well.
 For $\blamb\in \Pcal^{(d+1)\times (d-r)}$  let $\ep^\blamb:=\partorformep\in \bw^{r+1}A^*(\PP(F))$.
Again by~\cite{GatSant1}, the set $\{\ep^\blamb\,|\, \blamb\in\Pcal^{(r+1)\times (d-r)}\}$ is an
$A^*(X)$-basis of $A^*(G(r+1, F))$. Denote by $[\Omega_{\blamb}(\varrho_{r,d}^* F_\bullet)]$ the class
in $A^*(G(r+1, F))$ of the $F_\bullet$-Schubert variety $\Omega_{\blamb}(\varrho_{r,d}^* F_\bullet)$.

\bclm{\bf Theorem.} \label{kmpflak} {\em The following equality holds: \be [\Omega_{\blamb}(\varrho_{r,d}^*
F_\bullet)]=[\Omega_{(\lambda_r)}(\varrho_0^*F_\bullet)]\w[\Omega_{(1+\lambda_{r-1})}(\varrho_1^*
F_\bullet)]\w\ldots\w[\Omega_{( r+\lambda_0)}(\varrho_r^* F_\bullet)]=\ep^{\blamb}\label{eq:kmpflksv}
\ee
 modulo the identification of $A^*(G(r+1, F))$ with $\bw^{r+1}A^*(\PP(F))$.  \qed
} \eclm

\noindent
\begin{small}
Equality~(\ref{eq:kmpflksv}) is  an elegant and compact re-interpretation of the determinantal formula
of Schubert Calculus proven by Kempf and Laksov in \cite{K-L} to compute classes of degeneracy loci of
maps of vector bundles. This formula was first generalized in~\cite[(8.3)]{pragacz1}, see also~\cite
[Example 3.5 and Appendix 4]{Pragacz}. Then a far reaching generalization was obtained in~\cite{Fulton}
with help of the correspondences in flag bundles. In fact in~\cite{pragacz1}, the ${\mathcal P}$-ideals
of polynomials supported on degeneracy loci were studied, giving a deeper insight in enumerative
geometry of these loci.  Formula ~(\ref{eq:kmpflksv}) was basically discovered in~\cite{Gat1} for
trivial bundles. The present formulation is as in~\cite{NigroTesi}.
\end{small}

\smallskip
\noindent Let us sketch the proof of Theorem~\ref{kmpflak}. Set $\mu_j:=d-r-\lambda_{r-j}$, then
\[
\bmu:=(\mu_0,\mu_1,\ldots, \mu_r)\in\Pcal^{(r+1)\times (d-r)}.
\]
Denote $A_j:=N_{\mu_{r-j}+j+1}$ (see Section~\ref{defNh}), i.e $A_j$ fits into the exact sequence
\[
 0\sra A_j\sra F\sra F_{d-(j+\mu_{r-j})-1}\sra 0.
 \]
  Then
$0\subsetneq A_0\subsetneq A_1\subsetneq\ldots\subsetneq A_r$ is a flag of subbundles of $F_d$.
The Schubert variety
\[
\Omega(A_0,A_1,\ldots, A_r)=\{\Lambda\in G(r+1,F)\,|\, \Lambda\cap A_i\geq i\}
\]
coincides with $\Omega_{\blamb}(\varrho_{r,d}^*F_\bullet)$ defined
by~(\ref{eq:omlamb}), as a simple check shows. Formula 7.9 in~\cite{LakTh1}, which translates the
determinantal formula proven in~\cite{K-L}, implies
\[
[\Omega(A_0,A_1,\ldots, A_r)]=[\Omega(A_0)]\w[\Omega(A_1)]\w\ldots\w[\Omega(A_r)],
\]
which is thence equivalent to~(\ref{eq:kmpflksv}). \qed

\section{Wro\'nskians of Sections of  Grassmann  Bundles of Jets }\label{wsgb}

\claim{} The general framework of Section~\ref{wrsecbngen} shows that the notion of linear system can be
generalized into that  of pairs $(\gamma, F_\bullet)$, where $F_\bullet$ is a vector bundle on $X$
equipped with a filtration and $\gamma$ a section of the Grassmann bundle $G(r+1, F)$. This picture can
be fruitfully applied in the case of (families of) smooth complex projective curves of genus $g\geq 0$.
For the time being let $C$ be any one such, and let $ L\in Pic^d(C)$. In this section we shall denote by
$\varrho_d:J^dL\sra C$ the bundle of jets of $L\sra C$ up to the order $d$. Accordingly,
 for each $0\leq r\leq d$, we shall denote  $\varrho_{r,d}:G(r+1, J^dL)\sra C$  the Grassmann bundle
 of  $(r+1)$-dimensional
subspaces of fibers of $\varrho$. The natural filtration of $J^dL$ given by  the quotients $J^dL\sra
J^iL\sra 0$ , for $-1\leq i\leq d$, will be denoted $J^\bullet L$ (setting $J^{-1}L=0$).

\claim{}\label{osculating} The {\em kernel filtration} of $J^dL$ \be N_\bullet(L):  0\subset
N_1(L)\subset\ldots\subset N_d(L)\subset N_{d+1}(L)=J^dL\label{eq:kernelfiltration} \ee is defined
through the  exact sequence of vector bundles $ 0\sra N_h(L)\sra J^dL\sra J^{d-h}L\sra 0 $, where
$N_{h}(L)$ is a vector bundle of rank $h$. It will be also  called the {\em osculating flag} -- see
below and Section~\ref{wmapprojline}. The fiber of $N_h(L)$ at $P\in C$ will be denoted by $N_{h,P}(L)$.

As in the previous section, the $\blamb$-generalized  {\em Wro\'nskian} subvariety of $G(r+1, J^dL)$ is
$\Omega_{\blamb}(\varrho_{r,d}^*J^\bullet L)$, which has  codimension $|\blamb|$ in $G(r+1, J^dL)$.  By
virtue of Proposition~\ref{opensubset},   the space $\Gamma_{\tt triv}(\varrho_{r,d})$ of sections
$\gamma$ of $\varrho_{r,d}$ such that $\gamma^*\Scal_r$ is a trivial subbundle of $J^dL$, can be
identified with an open subset of $G(r+1, H^0(J^dL))$.  Hence $\gamma^*\Scal_r$ is of the form $C\times
U$ for some $U\in G(r+1, H^0(J^dL))$. As in section~\ref{wrsecbngen} we gain a Wro\'nski map:
\be
\Gamma_{\tt triv}(\varrho_{r,d})\lra \PP H^0(\wlbund),\label{eq:wrmpjbnd}
\ee
defined by $\gamma\mapsto W_0(\gamma)\, (\mod\,\CC^*)$. As we said, this map is the restriction to the
open subset   $\Gtrdb\subseteq G(r+1, H^0(J^dL))$  of the determinant map
\[
G(r+1, H^0(J^dL))\sra \PP H^0(\wlbund),
\]
sending $U$ to $t_r(u_0)\w t_r(u_1)\w \ldots \w t_r(u_r)\, (\mod\,\CC^*)$, where $(u_0,u_1,\ldots, u_r)$
is a basis of $U$ and $t_r$ denotes the epimorphism $J^dL\sra J^rL$.

\claim{}\label{sec53} We notice now that each $g^r_d(L)$, i.e.  $V\in G(r+1, H^0(L))$,  can be
seen in fact as an element of $\Gamma_{\tt triv}(\varrho_{r,d})$, because  $\Dcal_d:C\times V \sra J^dL$
realizes
$C\times V$ as a (trivial) vector subbundle of $J^dL$.  Indeed $D_dV:=\{D_dv\,|\, v\in V\}$ is an
$(r+1)$-dimensional subspace of $H^0(J^dL)$ because the map $J^dL\sra L\sra 0$ induces the surjection
$H^0(J^dL)\sra H^0(L)\sra 0$, see e.g.~\cite{CuGaNi}, and then $\Dcal_dv=0$ implies $v=0$.

We have thus an injective map $G(r+1, H^0(L))\hookrightarrow \Gamma_{\tt triv}(\varrho_{r,d})\subseteq G(r+1,
H^0(J^dL))$, sending $V$ to $\gamma_{D_dV}$, and
\[
W_0(\gamma_{\Dcal_dV}):=D_ru_0\w D_r u_1\w\ldots\w D_r u_r\,\,\mod\,\,\CC^*=W_0(V)
\]
which proves that our Wro\'nski map defined on $\Gamma_{\tt triv}(\varrho_{r,d})$,  which is in general
strictly larger than $G(r+1, H^0(L))$, coincides with the Wro\'nskian $W_0(V)$ defined in
section~\ref{wrseclbclass}. We are so in condition of defining {\em generalized Wro\'nskian subloci}.
Recall the  natural evaluation map
\[
\ev:C\times \Gamma_{\tt triv}(\varrho_{r,d})\lra G(r+1, J^dL)
\]
sending $(P, \gamma)\mapsto \gamma(P)$. If $\Omega_\blamb(\varrho_{r,d}^*J^\bullet L)$ is a generalized
Wro\'nski variety of $G(r+1, J^dL)$, then $\ev^{-1}(\Omega_\blamb(\varrho_{r,d}^*J^\bullet L))$ cuts the
locus of pairs $(P, \gamma)$ such that $\gamma(P)\in \Omega_\blamb(J^\bullet L)$. We also set
$\ev_P(\gamma)=\gamma(P)$,  for each $P\in C.$  It follows that the general section of any irreducible
component of $ \ev_P^{-1}(\Omega_\blamb(J^\bullet L)) $ is a section having $\blamb$ as a ramification
partition.

\claim{} The map $D_{d,P}:H^0(L)\sra J^d_PL$ sending $v\mapsto D_dv(P)$ is a vector space monomorphism.
If $V\in G(r+1, H^0(L))$, then $v\in V\cap D_{d,p}^{-1}(N_{h,P}(L))$ if and only if $D^hv(P)=0$, i.e.
if and only if $v$ vanishes at $P$ with multiplicity at least $h$. This explains the terminology
{\em osculating flag} used in Section~\ref{osculating}.

\claim{\bf Example.} \label{ex54} More details about  the present example  are in~\cite{GatSalehyan}.
Let $\gf$ be a proper
flat family of smooth projective curves of genus $g\geq 2$.  The {\em Hodge bundle} of
the family is  $\EE_\pi:=\pi_*K_\pi$. The vector bundle map over $\TS$
\begin{center}
$\pi^*\EE_\pi\sra J^{2g-2}K_\pi$
\end{center}
is injective and then it induces a section $\gamma_K:\TS\sra G(g, J^{2g-2}K_\pi)$. In this case the {\em
cuspidal locus} of $\gamma_K$,  which is by definition $\gamma_K^{-1}(\Omega_{1^{g-1}}(J^\bullet K))$,
coincides with the locus in $\TS$ of the Weierstrass points of the hyperelliptic fibers of $\pi$. With
the notation as in ~\ref{intthgrbund} and~\ref{kmpflak}, its class in $A^{g-1}(\TS)$ is given by
\[
[\gamma^{-1}_K(\Omega_{1^{g-1}}(J^\bullet K)]=\gamma_K^*[\Omega_{1^{g-1}}(J^\bullet K)]=
\gamma_K^*(\ep^0\w\ep^2\w\ldots\w\ep^g)
\]
and can be easily computed through a straightforward computation (see~\cite[Section~3]{GatSalehyan},
where the computation was performed for $g=4$). Since on each hyperelliptic fiber there are precisely
$2g+2$ Weierstrass points, the class of the hyperelliptic locus in $A^{g-2}(S)$ is given by
\[
[H]={1\over 2g+2}\cdot \pi_*\gamma_K^*(\ep^0\w\ep^2\w\ldots\w\ep^g)\,,
\]
which yields precisely the formula displayed  in~\cite[p.~314]{MumfordEnum}.

\claim{}  If $C=\PP^1$ and $L=O_{\PP^1}(d)$, then $\Gamma_{\tt triv}(\varrho_{r,d})$ coincides in this
case with $G(r+1, H^0(L))$ and our picture allows to rephrase in an elegant way the situation exposed in
the first part of \cite{EHcuspidal}. The {\em Wro\'nski map} $\Gamma_{\tt triv}(\varrho_{r,d})\sra \PP
H^0(O_{\PP^1}((r+1)(d-r)))$ coincides with~(\ref{eq:wronskimap}), modulo the identification of
$\Gamma_{\tt triv}(\varrho_{r,d})$ with $G(r+1, H^0(O_{\PP^1}(d)))$. In other words, when $C$ is not
rational, the theory exposed up to now  is a generalization of the  theory of linear systems on the
projective line, for which we want  to spend some additional words in a separate section.

\section{ Linear Systems on $\PP^1$  and the Intermediate Wro\'nskians}\label{wmapprojline}

In the case of linear systems $g^r_d$ defined on the projective line, the picture outlined in
Section~\ref{wsgb} gets simpler. However, even this case is particularly  rich of nice geometry
interacting with other parts of mathematics.

For the sake of brevity, denote by $L_d$ the invertible sheaf $O_{\PP^1}(d)$, i.e.  the unique line
bundle on $\PP^1$ of a fixed degree $d$. The elements of a basis $\bfx:=(x_0,x_1)$  of $H^0(L_1)$ can be
regarded  as homogeneous coordinates $(x_0:x_1)$ on $\PP^1$. Furthermore $H^0(L_d)=\Sym^dH^0(L_1)$, i.e.
$H^0(L_d)$ can be identified with the $\CC$-vector space generated by the monomials
$\{x_0^{i}x_1^{d-i}\}_{0\leq i\leq d}$, and a $g^r_d$ on $\PP^1$ is a point of $G(r+1, H^0(L_d))$. Any
basis $\bfv:=(v_0,v_1,\ldots, v_r)$ of $V\in G(r+1, H^0(L_d))$ defines a rational map
\be \varphi_V\,:\,\PP^1\rightarrow \PP^{r}\,, \quad P\mapsto (v_0(P):v_1(P):\ldots: v_r(P))\,.
\label{eq:projratcurv}\ee If $V$ has no  base points (that is, if $\dim V(-P)=\dim V-1$ for each
$P\in\PP^1$), then the image of~(\ref{eq:projratcurv}) is a non-degenerated (that is, not contained in
any hyperplane) rational curve of degree $d$ in $\PP^r$. In particular, if $r+1=\dim H^0(L_d)$, then
$V=H^0(L_d)$ and $\varphi_V(\PP^1)$ is nothing else than the rational normal curve of degree $d$. Each
curve  of degree $d$ in $\PP^r$  can be seen as the rational normal curve in $\PP H^0(L_d)$ composed
with a projection $\PP H^0(L_d) \rightarrow \PP^r$ whose center is a complementary linear subvariety of
$V\in G(r+1, H^0(L_d))$ (see e.g.~\cite{EHcuspidal}, \cite{KhSo}).

Keeping the notation of Section~\ref{wsgb}, let $\varrho_d:J^dL_d\sra \PP^1$ be the bundle of $d$-jets
of $L_d$. Then $\Dcal_d:\PP^1\times H^0(L_d)\sra J^dL $   (cf.~(\ref{eq:vbmsr}) is an injective morphism
between vector bundles of the same rank, that is, an isomorphism.  In particular, the map
 \be
 \left\{\matrix{\Dcal_{d,P}&:&H^0(L_d)&\lra& J^d_PL_d\cr
 {}&{}&P&\longmapsto&D_dv(P)}\right.
\ee
is an isomorphism of  vector spaces, for each $P\in\PP^1$. We define the {\em osculating flag} at $P$ of
$H^0(L_d)$,
 \[
\Fcal_{\bullet, P}: 0\subset \Fcal_{1,P}\subset \ldots\subset \Fcal_{d,P}\subset \Fcal_{d+1,P}=J_P^dL,
 \]
by setting  (cf. \ref{osculating})
 \[
 \Fcal_{h,P}=\Dcal_{d ,p}^{-1}(N_{h,P}(L))\subseteq H^0(L_d).
 \]
 In other words, $v\in V\cap \Fcal_{h,P}$  if and only
 if $v$ vanishes at $P$ with multiplicity at least $h$, that is,  $\Dcal_hv(P)=0$.
 In fact, $\Fcal_{h,P}$ may be identified with the vector subspace  of  the
 homogeneous polynomials of $H^0(L_d)$  that vanish at $P$ with multiplicity at least $h$. Yet another
 interpretation of  $\Fcal_{h,P}$ is
 the set of all hyperplanes of $\PP H^0(L_d)$ intersecting the rational normal curve
 in $\PP H^0(L_d))$ at $P$ with multiplicity at least $d-h$.

\claim{} \label{61} The Riemann-Roch formula shows that $h^0(L_d)=h^0(J^dL_d)$; thus the injective
``derivative map" $D_d: H^0(L_d)\sra H^0(J^dL_d)$ is an isomorphism which itself induces  a
biholomorphism:
\[
G(r+1, H^0(L_d))\sra G(r+1, H^0(J^dL_d)).
\]
So one concludes that $\Gamma_{\tt triv}(\varrho_{r,d})=G(r+1, H^0(J^dL_d))\cong G(r+1, H^0(L_d))$
parameterizes all the $g^r_d$'s on $\PP^1$ (with base points or not). In particular it is compact.

For $V\in G(r+1, H^0(L_d))$, denote by $\gamma_V$ the corresponding element of $\Gamma_{\tt
triv}(\varrho_{r,d})$. The evaluation morphism  $\PP^1\times G(r+1, H^0(L_d))\sra G(r+1, J^dL_d)$ maps
$(P,V)$ to $\gamma_V(P)\in G(r+1,J^d_PL_d)$.

\claim{} By~\ref{61},  the {\em Wro\'nski map}  $\gamma\mapsto W_0(\gamma)$ (see (\ref{eq:wrmpjbnd}))
coincides with the Wro\'nski map~(\ref{eq:wronskimap}) of Section~\ref{S39}:
 \be
G(r+1, H^0(L_d))\sra \PP H^0(L_{(r+1)(d-r)})\,,\quad V\mapsto W_0(V).\label{eq:wrmpp1}
\ee
It is a finite surjective morphism (see e.g.~\cite{EHcuspidal}, \cite{KhSo}, \cite{Scherb1}). Its degree
$N_{r,d}$ is precisely the Pl\"ucker degree of the Grassmannian $G(r+1, d+1)$:
\[
N_{r,d}=\int \sigma_{(1)}^{(r+1)(d-r)}\cap[G(r+1, d+1)]={1!2!\ldots r!\cdot (r+1)(d-r)!\over
(d-r)!(d-r+1)!\cdot\ldots \cdot d!}\,.
\]
Thus, given a homogeneous polynomial $W$ of degree $(d-r)(r+1)$ in two indeterminates $(x_0,x_1)$, there
are at most $N_{r,d}$  distinct $g^r_d$'s having $W$ as a Wro\'nskian. The number $N_{r,d}$  was
calculated by Schubert himself in 1886, cf.~ \cite{schubert} and~\cite[p.~274]{Fu1}. In the case of real
rational curves, the degree of the Wro\'nski map was obtained by L.~Goldberg for $r=1$  (\cite{Gold}),
and for any $r\geq 1$ by A.~Eremenko and A.~Gabrielov (\cite{EreGab1}). For more considerations on real
Wro\'nski map see also~\cite{KhSo}.

\claim{}\label{inters1}    For any partition $\blamb\in \Pcal^{(r+1)\times(d-r)}$ define
 \[
 \Omega_\blamb(P):=\Omega_\blamb(\Fcal_{\bullet,P})\subseteq G(r+1, H^0(L_d)).
 \]
 It is a Schubert variety of codimension $|\blamb|$ in $G(r+1, H^0(L_d))$. If $\blamb(V,P)$ is the order
 partition of $V$ at $P$  (see Section~\ref{orderpartition}) then
\[
V\in \Omega^{\circ}_{\blamb(V,P)}(P)\subseteq \Omega_{\blamb(V,P)}(P)\,,
\]
 and  $P$ is a $V$-ramification point if and only if  $|\blamb(V,P)|>0$.
 The Wro\'nskian $W_0(V)$ of $V$ vanishes exactly at the $V$-ramification points. The  total weight of
 the $V$-ramification points equals the dimension of  $G(r+1, H^0(L_d))$ (one can see that by putting $g=0$ in
 ~(\ref{eq:brillsegre})).

 Let $\{(\underline{P},\underline{w})\}:=\{(P_0,w_0),( P_1,w_1), \ldots, (P_k,w_k))$ be a $k+1$-tuple
 of pairs where  $P_i\in \PP^1$  and $w_i$'s are positive integers such that
 \be
 \sum_{i=1}^k w_i=(r+1)(d-r).\label{eq:sumwts}
 \ee
 Thus, in notation of Section~\ref{intrwrsec}, if
 \[
 D_{w_i-1}W_0(V)\in H^0(J^{w_i-1}L_{(r+1)(d-r)})
 \]
 vanishes at $P_i$, for every $0\leq i\leq k$, then  $P_0,P_1,\ldots, P_k$ are exactly the
 ramification points of $V$,  each one of weight $wt_V(P_i)=w_i=|\lambda(P_i, V)|$. We have
 \begin{eqnarray}
 V&\in&\Omega^\circ_{\blamb(V,P_0)}(P_0)\cap \Omega^\circ_{\blamb(V,P_1)}(P_1)\cap\ldots\cap
 \Omega^\circ_{\blamb(V,P_k)}(P_k)=\nonumber\\ &=&\Omega_{\blamb(V,P_0)}(P_0)
 \cap\Omega_{\blamb(V,P_1)}(P_1)\cap
 \ldots\cap \Omega_{\blamb(V,P_k)}(P_k)\,.\label{eq:intramprt}
 \end{eqnarray}

Condition~(\ref{eq:sumwts}) means that  the "expected dimension" of the intersection
(\ref{eq:intramprt}) is zero. Intersections of Schubert varieties associated with the osculating flags
of the normal rational curve were first studied by D.~Eisenbud and J.~Harris in the eighties,
\cite{EHcuspidal}. In particular, they showed that the intersection (\ref{eq:intramprt}) is
zero-dimensional indeed, and hence the number of  distinct elements in the intersection is at most
\[
\int_{G(r+1, H^0(L_d))}\sigma_{\blamb(P_0,V)}\cdot\sigma_{\blamb(P_1,V)}\cdot \ldots
\cdot\sigma_{\blamb(P_k,V)}\cap [G(r+1, H^0(L_d))]\,,
\]
where  $\sigma_{\blamb}$ is the Schubert cycle  defined by the equality $\sigma_{\blamb}\cap[G(r+1,
H^0(L_d))]=[\Omega_\blamb]$. This fact was used in~\cite{NBFSC} to deduce explicit formulas (and a list
up to $n=40$) for the number of space rational curves of degree $n-3$ having $2n$ {\em hyperstalls} at
$2n$ prescribed points.

\claim{\bf Preimages of the Wro\'nski Map.} \label{pwmiw} Notice that if  $P\in \PP^1$ is a base point
of $V$, it occurs in the  $V$-ramification locus as well, and  the Wro\'nskian vanishes at it with
weight $(r+1)$. The set $B_P$ of linear systems having $P$ as base point is a closed subset of $G(r+1,
H^0(L_d))$ of
 codimension $(r+1)$. In fact $B_P:=\ev_P^{-1}(\Bcal(\varrho_{r,d}^*J^\bullet L_d))$, which is a closed
 subset of codimension $(r+1)$ (cf. Section~\ref{bsptlc}).

Let $\{(\underline{P},\underline{w})\}$ be as in~\ref{inters1}. Denote by $\Gcal_{r,d}(\underline{P})$
the set of all $V\in G(r+1, H^0(L_d))$ whose base locus contains no $P_i$, $0\leq i\leq r$. It is an
open dense subset of codimension $(r+1)$,
\[
\Gcal_{r,d}(\underline{P})=G(r+1, H^0(L_d))\setminus( B_{P_0}\cup B_{P_1}\cup\ldots\cup B_{P_k})\,.
\]
Consider now a $(k+1)$-tuple of partitions
\[
\vec{\blamb}=(\blamb_0,\blamb_1,\ldots,\blamb_k)\,, \ \ |\blamb_j|=w_j, \ \ 0\leq j\leq k.
\]
We shall write:
\[
\blamb_j:=\lambda_{j,0}\geq \lambda_{j,1}\geq\ldots\geq \lambda_{j,r}\,.
\]
The elements of
\be
I(\vec{\blamb}, \underline{P})=\Omega_{\blamb_0}(P_0)\cap\Omega_{\blamb_1}(P_1) \cap\ldots\cap
\Omega_{\blamb_k}(P_k)\cap \Gcal_{r,d}(\underline{P})\,\subset G(r+1, H^0(L_d)) \label{eq:intersection}
\ee
correspond to the base point free linear systems ramifying at $\underline{P}$ according to
$\vec{\blamb}$.

The problem of determining $I(\vec{\blamb}, \underline{P})$ leads to interesting analytic considerations
related with Wro\'nskians. Up to a projective change of coordinates, it  is not restrictive to assume
that $P_0=\infty:=(0:1)$. Using the coordinate $x=x_1/x_0$, the osculating flag at $\infty$   shall be
denoted by $\Fcal_{\bullet,\infty}$. Accordingly,   the partition $\blamb_0$ will be renamed
$\blamb_\infty$. Notice that $\Fcal_{j,\infty}$ coincides with the vector space $\Poly_j$ of the
polynomials of degree at most $j$ in the variable $x$: in fact  a polynomial $P(x)$ (thought of as the
affine representation of a homogeneous polynomial of degree $d$ in two variables)   vanishes at $\infty$
with multiplicity $j$ if and only if it has degree $d-j$. For $V\in I(\vec{\blamb}, \underline{P})$, let
\[
W_V(x):={W_0(V)\over x_0^{(r+1)(d-r)}}
\]
be the representation of the $W_0(V)$ in the affine open subset of $\PP^1$ defined by $x_0\neq 0$. The
degree of the polynomial $W_V(x)$ is less or equal than $(r+1)(d-r)$, because of possible ramifications
of $V$ at $\infty$. We have
\be
W_V(x)=(x-z_1)^{w_1}\cdot\ldots\cdot (x-z_k)^{w_k}, \label{eq:Wr}
\ee
where $z_i:=x(P_i)$ are the values of the coordinate $x$ at $P_i\in \PP^1$; $\sum_{i=1}^k w_i=\deg
W_V(x)\leq (r+1)(d-r)$.

For  a basis $\bfv=(v_0,v_1,\ldots, v_r)$ of $V$, consider $f_i:={v_i/x_0^d}$ and write
$\bff=(f_0,f_1,\ldots, f_r)$. According to~(\ref{eq:wrrpuform}), one writes $ W_V(x)=\bff\w
D\bff\w\ldots \w D^r\bff $, where
\[
D^j\bff= \left({d^jf_i\over dx^j}\right)_{0\leq i\leq r}.
\]
The space $V$ can be realized as the solution space of the following differential equation
\be
E_V(g)=\left|\matrix{g&f_0&f_1&\ldots&f_r\cr Dg&Df_0&Df_1&\ldots&Df_r\cr
\vdots&\vdots&\vdots&\ddots&\vdots\cr D^rg&D^rf_0&D^rf_1&\ldots&D^rf_r \cr
D^{r+1}g&D^{r+1}f_0&D^{r+1}f_1&\ldots&D^{r+1}f_r}\right|=0\,. \label{eq:EV}
\ee

\claim{\bf Intermediate Wro\'nskians.}\label{intwr} For  $V\in I(\vec{\blamb},\underline{P})$, denote by
$V_\bullet$  the  flag obtained by the intersection of $V$ and  $\F_{\bullet, \infty}$:
\begin{equation}\label{V}
V_\bullet\,=\,\left\{ V_0\subset V_1\subset V_2 \subset \dots \subset V_r=V\,\right\},\quad \dim
V_j=j+1,
\end{equation}
all the polynomials in $V_j$ have degree $\leq d_j$, where $0\leq d_0< d_1< \dots < d_r\leq d$ is the
order sequence of $V$ at $P_0$ (cf. Section~\ref{orderpartition}). Recall that $V$ has no base point and
$W_V(x)$ as in (\ref{eq:Wr}).

Define {\em the $j$-th intermediate Wro\'nskian of} $V$ as $W_j(x):=W_{V_j}(x)$, the Wro\'nskian of
$V_j$, $0\leq j\leq r$. In particular, the $r$-th intermediate Wro\'nskian coincides with $W_V(x)$.
Non-vanishing properties of intermediate Wro\'nskians have been recently investigated in an analytic
context in~\cite{camporesi} and~\cite{campodisc} to study factorizations of linear differential
operators with non-constant $\CC$-valued coefficients.

Intermediate Wro\'nskians are important because every $V\in G(r+1, {\Poly}_d)$ is completely determined
by the set of its intermediate Wro\'nskians $W_0(x),\dots\,, W_r(x)$. Indeed, the ODE (\ref{eq:EV}) can
be rewritten as follows:
\[
\frac{d}{dx}\,\frac{W^2_{r}(x)}{W_{r-1}(x)W_{r+1}(x)}\cdot  \dots \cdot
\frac{d}{dx}\,\frac{W^2_2(x)}{W_3(x)W_1(x)}\,\frac{d}{dx}\cdot\frac{W^2_1(x)}{W_2(x)W_0(x)}\cdot
\frac{d}{dx}\,\frac{g(x)}{W_1(x)}=0\,.
\]
By~\cite[Part VII, Section~5, Problem 62]{PSz},  one can
take as a basis of $V$ the following set of $r+1$ linearly independent solutions of~(\ref{eq:EV}):
\begin{eqnarray*}
 g_0(x) & = & W_0(x)\,, \\
 g_1(x)&=&W_0(x)\int^x\frac{W_0W_2}{W_1^2}\,,\\
 g_2(x)&=&W_0(x)\int^x\left(\frac{W_0(\xi)W_2(\xi)}{W_1^2(\xi)}
 \int^{\xi}\frac{W_1W_3}{W_2^2}\right)\,,\\
 \dots & \ & \dots \quad \ \ \quad \dots\quad \ \ \quad \dots \quad \ \ \quad \dots  \\
 g_r(x) & = &
W_0(x)\int^x\left(\frac{W_0(\xi)W_2(\xi)}{W_1^2(\xi)}\int^{\xi}\left(
\frac{W_1(\tau)W_3(\tau)}{W_2^2(\tau)}\int^{\tau}\dots\int^{\eta}\left.
\frac{W_{r-1}W_{r+1}}{W^2_{r}}\right.\right)\dots\right)\,.
\end{eqnarray*}

Define now   polynomials $Z_0(x), Z_1(x),\ldots, Z_r(x)$ through the following formula:
\be\label{Z} Z_i(x)=\prod_{j=1}^k(x-z_j)^{m_j(i)},\quad 0\leq i\leq r\,,
\ee
where
$$
m_j(i)=\lambda_{j, r}+\lambda_{j,r-1}+\ldots+\lambda_{j, r-i},\qquad 1\leq j\leq k.
$$
In particular  $Z_r(x)=W_V(x)$.

\bclm{\bf Lemma}\label{l1-3} (\cite{Scherb2}) {\em The ratio $T_{r-i}(x):=W_i(x)/Z_i(x)$
is a polynomial of degree}
\be\label{k}
(i+1)(d-i)-\sum_{l=0}^{i}\lambda_{r-l,\infty}- \sum_{j=1}^{k} m_j(i).
\ee
\eclm
In particular, $T_0(x)=1$. Thus we have $W_{r-j}(x)=T_j(x)Z_{r-j}(x)$, \ $0\leq j\leq r$. The roots of
$T_j(x)$ are said to be {\it the additional roots} of the $(r-j)$-th  intermediate Wro\'nskian. If
(\ref{eq:intersection})  contains more than one element, then the intermediate Wro\'nskians of these
elements all differ by the additional roots.

\claim{\bf Non-degenerate planes.} (\cite{Scherb2}) The intersection (\ref{eq:intersection}) contains
some distinguished elements,  called {\it non-degenerate planes}. Denote by $\Delta(f)$ the discriminant
of a polynomial $f(x)$ and by $\Res(f,g)$ the resultant of polynomials $f(x),\, g(x)$.

\medskip\noindent{\bf Definition.}  {\em We call $V\in I(\vec{\blamb}, \underline{P})$
a {\sl non-degenerate plane}
if the polynomials $T_0(x),\dots\,, T_{r-1}(x)$
\begin{itemize}
\item[i)]
do not vanish at the ramification points $P_1,\ldots, P_k$, i.e. $T_i(z_j)\neq 0$ for all
$0\leq i\leq r-1$ and all $1\leq j\leq k$;
\item[ii)]
do not have multiple roots: $\Delta (T_i)\,\neq\, 0$,  for all $0\leq i\leq r$;
\item[iii)]
For each $1\leq i\leq r$ , $T_i$  and $T_{i-1}$  have no common roots:  $\Res(T_i,T_{i-1})\,\neq\,0$.
\end{itemize}
}
\claim{\bf Relative discriminants and resultants.}\ Non-degenerate planes correspond to critical
points of a certain {\em generating function} which can be described  in terms of {\em relative}
discriminants and resultants.
 For fixed $\bfz=(z_1,\dots,z_k)$,  any monic polynomial $f(x)$
can be written in a unique way as the product of two monic polynomials $T(x)$ and $Z(x)$ satisfying
\be\label{TZ}
f(x)=T(x)Z(x), \ \   T(z_j)\neq 0,\ \  Z(x)\neq 0 {\rm\ \ for\ any\ } x\neq z_j, \ \ 1\leq j\leq k.
\ee
One defines  {\em  the relative discriminant of $f(x)$  with respect to $z$} as being
\[
\Delta_{\bfz}(f)=\frac{\Delta(f)}{\Delta(Z)}=\Delta(T){\rm (Res}(Z,T))^2,
\]
and {\em the relative resultant   of  $f_i(x)=T_i(x)Z_i(x)$, $i=1,2$, with respect to $\bfz$} as
$${\rm Res}_{\bfz}(f_1, f_2)=\frac{{\rm Res}(f_1,f_2)}{{\rm Res}(Z_1,Z_2)}=
{\rm Res}(T_1, T_2){\rm Res}(T_1, Z_2){\rm Res}(T_2, Z_1).$$

\medskip
If $V$ is a non-degenerate plane in  $I(\vec{\blamb}, {\underline{P}})$ given by
(\ref{eq:intersection}), then the decomposition $W_i(x)=T_{r-i}(x)Z_i(x)$ is exactly the same as
displayed in (\ref{TZ}).
 The {\em generating function} of $I(\vec{\blamb}, {\underline{P}})$
is a rational function such that its critical points determine the non-degenerate elements in such an
intersection.  Its expression is  (see \cite{Scherb2} ):

\begin{equation}\label{eq:Phi-w}
\Phi_{ (\vec{\blamb}, \bfz)} (T_1,\dots, T_{p-1})= \frac {\Delta_\bfz(W_0)\cdot \ldots\cdot
\Delta_\bfz(W_{r-1})} {{\Res}_\bfz(W_1,W_2)\cdot\ldots\cdot {\Res}_\bfz(W_{r-1},W_r)}
\end{equation}

\noindent Part of the following theorem was originally obtained by A.~Gabrielov (unpublished), along his
investigations of the Wro\'nski map.

\bclm{\bf Theorem}\label{T}  ( \cite{Scherb2}) {\em There is a one-to-one correspondence between the
critical points with non-zero critical values of the function $\Phi_{ (\vec{\blamb}, \bfz)} (T_0,\dots,
T_{r-1})$ and the non-degenerate planes in  the intersection $I(\vec{\blamb}, {\underline{P}})$ given by {
{\em (\ref{eq:intersection})}}.
}
\eclm

Namely, every such critical point defines  the intermediate Wro\'nskians, and hence a non-degenerate
plane, see~\ref{intwr}. Conversely, for every non-degenerate plane one can calculate the intermediate
Wro\'nskians, and the corresponding polynomials $T_i(x)$ supply a critical point with a non-zero
critical value  of the {\em generating function} (\ref{eq:Phi-w}).

\claim{\bf Relation to Bethe vectors in the Gaudin model} (see \cite{ MV, Scherb, Scherb2}). Once one
re-writes (\ref{eq:Phi-w}) in terms of unknown roots of the polynomials $T_j$'s, the generating function
turns into the {\em master  function}  associated with the {\em Gaudin model} of statistical mechanics.

In the Gaudin model, the partitions $\blamb_j$,  $1\leq j\leq k$, of  Section~\ref{pwmiw} are the
highest weights of ${\frak sl}_{r+1}$-representations, and the $j$-th representation is marked by the
point $P_j$. Recall that $\blamb_\infty$ is the partition related to  $P_0:=(0:1)\in \PP^1$, after
renaming $\blamb_0$, see Section~\ref{pwmiw}. Denote by $\blamb_\infty^*$ the partition dual to
$\blamb_\infty$.  Certain commuting linear operators, called {\em Gaudin Hamiltonians}, act in the
subspace of singular vectors of the weight $\blamb_\infty^*$  in the tensor product of the ${\frak
sl}_{r+1}$-representations of the weights $\blamb_1,\ldots, \blamb_k$, and one looks for a common
eigenbasis of the Gaudin Hamiltonians.

  The {\em Bethe Ansatz} is a method to look for common eigenvectors.  It gives a family of vectors
  of the required weight $\blamb_\infty^*$  meromorphically
 depending on a number of auxiliary complex parameters. The {\em Bethe system} is a system of
 equations on these parameters, and any member of the family that corresponds to a solution of the Bethe
 system is a common singular eigenvector of the Gaudin Hamiltonians called the {\em Bethe vector}.

It turns out that the Bethe system coincides with the system on critical points with non-zero critical value
 of the function $\Phi_{\vec{\blamb},\bfz,}$. In other words, the auxiliary complex parameters are exactly
 the  additional roots of the intermediate Wro\'nskians! Thus every non-degenerate plane of
 (\ref{eq:intersection}) defines a Bethe vector and vice versa.

This link has led to an essential progress in studies of the Gauidin model as well as in algebraic
geometry (e.g., Shapiro-Shapiro conjecture), see \cite{MTV} and references therein.

\section{Wronskians of   (hyper)elliptic involutions}

The main reference for this section is~\cite{CuGaNi}.
\claim{} Let $C$ be a smooth projective curve of genus $g\geq 1$ and let  $\Mcal\in Pic^2(C)$ such that $h^0(C,\Mcal)=2$. Then $C$ is elliptic if $g=1$, and $\Mcal=O_C(2P)$ for some $P\in C$, or  hyperelliptic if $g\geq 2$ and then $\Mcal$ is of the form $O_C(2P)$,  where $P$ is some of the Weierstrass points of $C$. If $K$ is the canonical bundle of a hyperelliptic curve, then  $K=\Mcal^{\otimes g-1}$ ~\cite[I-D9, p.~41]{ACGH}.  (Hyper)elliptic curves are, in a sense, the closest example to rational curves, as they are   double ramified coverings of the projective line. How does the wronski map look like in this case?
One of the main result of~\cite{CuGaNi} is that the  {\em extended Wronski map}
\[
\Gamma_{\tt t}(\rho_{1,2})\sra \PP H^0(C, \Mcal^{\otimes 2}\otimes K)= \PP H^0(C, \Mcal^{\otimes g+1})
\]
is dominant (notation as in~Section~\ref{osculating}).  The proof relies on producing an explicit basis of the space $H^0(\Mcal^{\otimes g+1})$.
One also shows that if ${\mathbf v}=(v_0,v_1)$ generates $H^0(\Mcal)$, then a basis of
$H^0(\Mcal^{\otimes g+1})$ is formed by
\[
\{v_0^{g+1-i}v_1^{i}, W_0(\mathbf v)\}_{0\leq i\leq  g+1}
\]
where $W_0({\mathbf v})$ is precisely the Wronskian of the chosen basis ${\mathbf v}$ of $H^0(\Mcal)$. More generally one has:
\bclm{\bf Theorem.} (Cf.~\cite[6.6]{CuGaNi}) {\em  The following direct sum decomposition holds:
\be
H^0(\Mcal^{\otimes a})=\Sym^{a-g-1}H^0(\Mcal)\cdot W({\bm\lambda})\oplus \Sym^aH^0(\Mcal),
\label{eq:mainresult}
\ee
where 
$
\Sym^jH^0(\Mcal)\cdot W({\bm\lambda})
$
is the image of $\Sym^jH^0(\Mcal)$ in $H^0(\Mcal^{\otimes g+1+j})$ through the multi- \linebreak plication--by--$W({\bm\lambda})$ map
$
H^0(\Mcal^{\otimes j})\sra H^0(\Mcal^{\otimes g+1+j}).
$
}
\eclm

\bclm{\bf Example.}\label{elptcex} Let $C$ be an elliptic curve and let $L=O_C(2P_0)$ for some point $P_0\in C$. Let $P_1,P_2,P_3$ be the remaining ramification points of $H^0(L)$. By definition of ramification points  there exists $v_0$ and $v_1$ vanishing at $P_0$ and $P_1$ with multiplicity $2$. Then  $\bfv:=(v_0,v_1)$ is a basis of $H^0(L)$ and  the map
\begin{center}
$(W_0(\bfv):v_0^2:v_0v_1:v_1^2):C\lra \PP^3$
\end{center}
realizes the elliptic curve as quartic curve in $\PP^3$,  a complete intersection of two quadrics. 

If $(X_0:X_1:X_2:X_3)$ are homogeneous coordinates of $\PP^3$, one of the two quadrics is the cone in $\PP^3$ of equation $X_1X_3-X_2^2=0$. To find the second quadric one argues as follows. Let $v_2$ and $v_3$ such that $Dv_2(P_2)=0$ and $Dv_3(P_3)=0$. Since  $(v_0,v_1)$ form a basis of $H^0(L)$, there are $a,b\in\CC^*$ such that  $v_2=v_1-av_0$ and $v_3=v_1-bv_0$.

The product $v_0v_1v_2v_3$ is a section of $L^{\otimes 4}$ which vanishes at each ramification point with multiplicity $2$. Since the Wronskian $W_0(\bfv)$ vanishes at each $P_i$ with multiplicity $1$, up to a multiplicative constant one has the the relation
\[
W(\bfv)^2=v_0v_1(v_1-av_0)(v_1-bv_0)=v_0v_1(v_1^2-(a+b)v_1v_0+abv_0^2)
\]
which shows that the image of $C$ is contained in the quadric
\[
X_0^2-X_2X_3-(a+b)X_2^2+abX_1X_2=0
\]
Again, one may observe that  defining $x=v_1/v_0$ and $y=W(\bfv)/v_0^2$
one gets the equation:
\[
y^2=x(x-a)(x-b)
\]
which is the classical affine {\em Weierstrass  equation} for an elliptic curve (which is then also a {\em Wronski equation}, too). One then sees  that its  natural
compactification lives  in the weighted projective space $\PP(2,1,1)$ or as  a quadric section of a quadric cone of $\PP^3$. Notice that on a trivializing set $U$ of $C$,  with local parameter $x$,  one may write  $x=f_1/f_0$, where $f_0,f_1$ are local holomorphic functions representing $v_0,v_1$ on $U$, respectively. Then one sees that:
\[
(x_{|_U})':={d\over dz}\left(\left.{v_1\over v_0}\right|_{{U}}\right)={d\over dz}\left({f_1\over f_0}\right)={W(f_0,f_1)\over f_0^2}=\left.{W(\bfv)\over v_0^2}\right|=y_{|_U}
\]
i.e. $y=x'$ which is of course compatible with the fact that the ``parametric'' equations of an elliptic curve in the affine plane are given by $x=\wp_\Lambda(z)$ and $y=\wp_\Lambda'(z)$, where $\wp_\Lambda$ is the Weierstrass $\wp$-function associated to some lattice of  $\CC$.
\eclm
\claim{} In~\cite{CuGaNi}  one shows as in Example~\ref{elptcex} that a hyperelliptic curve of genus $g\geq 2$ satisfies a Weierstrass type equation  expressing the relation between the Wronskian and the product of the sections vanishing twice at its Weierstrass points. It is also intrinsically shown, again using the notion of Wronskian,  that the hyperelliptic curve of genus $g\geq 2$ can be 
realized, for each $a\geq 0$, 
as a curve of degree $2(g + 1 + a)$ lying on a rational normal scroll $S(a, g + 1 + a)$ of $\PP^{g+2+a}$. If $a=0$ it is a quadric section of a cone of degree $g+1$ in $\PP^{g+2}$ -- see also~\cite{transcanonical}.

\section{Linear ODEs and Wro\'nski--Schubert Calculus}\label{finalsec}
This last section surveys and announces the results of~\cite{GatScherb1}, an attempt  to reconcile the
first part of this survey, regarding Wro\'nskians of fundamental systems of solutions of linear ODEs,
with the geometry described  in the last four sections. The main observation is that Schubert cycles of
a Grassmann bundle  can be described through Wro\'nskians associated with a fundamental system of
solutions of a linear ODE. \claim{}  Let us work in the category of  (not necessarily finitely
generated) associative commutative $\QQ$-algebras with unit. Let $A$ be such a $\QQ$-algebra. We denote
by $A[T]$ and $A[[t]]$  the corresponding $A$-algebras of polynomials and of formal power series,
respectively (here $t$ and $T$ are indeterminates over $A$). For $\phi=\sum_{n\geq 0}a_nt^n\in A[[t]]$,
we write $\phi(0)$ for the ``constant term'' $a_0$. If $P(T)\in A[T]$ is a polynomial of degree $r+1$,
we denote by $(-1)^ie_i(P)$ the coefficient of $T^{r+1-i}$,  for each $0\leq i\leq r+1$; for instance,
if  $P$ is monic, $e_0(P)=1$, we have:
\[
P(T)=T^{r+1}-e_1(P)T^r+\ldots+(-1)^{r+1}e_{r+1}(P).
\]

\noindent Let $B$ be another $\QQ$-algebra. Each $\psi\in Hom_\QQ(A,B)$ induces two obvious
$\QQ$-algebra homomorphisms, $A[T]\sra B[T]$ and $A[[t]]\sra B[[t]]$, the both are also denoted by
$\psi$. The former is defined by $e_i(\psi(P))=\psi(e_i(P))$ and the latter by $\sum_{n\geq 0}
a_nt^n\mapsto \sum_{n\geq 0}\psi(a_n)t^n$.

\claim{} Let $E_r:=\QQ[e_1,e_2,\ldots, e_{r+1}]$ be the polynomial $\QQ$-algebra in the set of  indeterminates
$(e_1,\ldots, e_{r+1})$. We call
\[
U_{r+1}(T)=T^{r+1}-e_1T^r+\ldots+(-1)^{r+1}e_{r+1}
\]
the {\em universal monic polynomial} of degree $r+1$. Thus $e_i(U_{r+1}(T))=e_i$ for all $0\leq i\leq
r+1$.

Let $\bar\bfh:=(h_0,h_1,h_2,\ldots, h_r, h_{r+1}, \ldots)$ be the sequence in $E_r$ defined by the
equality of formal power series:
\begin{eqnarray*}
\sum_{n\geq 0}h_nt^n={1\over 1-e_1t+\ldots+(-1)^{r+1}t^{r+1}}=1+\sum_{n\geq
1}(e_1t-e_2t^2+\ldots+(-1)^re_{r+1}t^{r+1})^n.
\end{eqnarray*}
One gets $h_0=1$, $h_1=e_1$,  $h_2=e_1^2-e_2$, \ldots\, . In general $h_n=\det(e_{j-i+1})_{1\leq i,j\leq
n}$ (see~\cite[p.~264]{Fu1}).

For any $(r+1)$-tuple or a sequence $\bar\bfa=(a_0, a_1, \ldots)$ of elements of any $E_r$-module, we
set
\be
U_0(\bfa)=a_0\,, \quad U_i(\bfa)=a_i-e_1a_{i-1}+\ldots+(-1)^ie_ia_0\,, \  1\leq i\leq r.\label{eq:Ua}
\ee
 Although only $a_0,a_1,\ldots, a_r$ appear in (\ref{eq:Ua}), we prefer  to define
$U_j$ also for sequences. We have $U_i(\bar\bfh)=0$ for all $1\leq i\leq r$.

 \claim{} Let
$\bfx:=(x_0,x_1,\ldots, x_r)$ and $\bar\bff:=(f_n)_{n\geq 0}$ be two sets of indeterminates over $\QQ$.
Let
\[
E_r[\bfx,\bar\bff]:=E_r[x_0,x_1,\ldots, x_r; f_0,f_1,\ldots]
\]
be the $\QQ$-polynomial algebra and $E_r[\bfx,\bar\bff][[t]]$ the corresponding algebra   of formal
power series. Denote by $D:=d/dt$ the  usual formal derivative of formal power series. Its $j$-th
iterated is:
\[
D^j\left(\sum_{n\geq 0}a_n{t^n\over n!}\right)=\sum_{n\geq 0}a_{n+j}{t^{n}\over n!}\,, \ \ a_m\in
E_r[\bfx,\bar\bff].
\]
Evaluating the polynomial $U_{r+1}$ at $D$ we get the {\em universal differential operator:}
\[
U_{r+1}(D)=D^{r+1}-e_1D^r+\ldots+(-1)^{r+1}e_{r+1}.
\]
Let $f:=\sum_{n\geq 0}f_n{t^n\over n!}\in \QQ[\bar\bff][[t]]\subseteq E_r[\bfx,\bar\bff][[t]]$. Consider
the universal Cauchy problem for a linear ODE with constant coefficients: \be
\left\{\matrix{U_{r+1}(D)y&=&f,&\cr\cr D^iy(0)&=&x_i,&\,\,\,0\leq i\leq r.} \right.\label{eq:cauchyprob}
\ee
We look for  solution of~(\ref{eq:cauchyprob}) in  $E_r[\bfx,\bar\bff][[t]]$.

\bclm{\bf Theorem.} \label{nhunsol} {\rm (\cite{GatScherb1})} \ {\em Let $\sum_{n\geq 0}p_n\cdot t^n\in
E_r[\bfx,\bar\bff][[t]]$ be defined by: \be \sum_{n\geq
0}p_nt^n={U_0(\bfx)+U_1(\bfx)t+\ldots+U_r(\bfx)t^r+\sum_{n\geq r+1}f_{n-r-1}t^n\over
1-e_1t+\ldots+(-1)^{r+1}e_{r+1}t^{r+1}}\,,\label{eq:genfunct} \ee where $U_j$ are as in (\ref{eq:Ua}).
Then
\be
g:=\sum_{n\geq 0}p_n{t^n\over n!}\label{eq:univsol}
\ee
is the unique solution of the  Cauchy problem~(\ref{eq:cauchyprob}). \qed } \eclm

The {\em universality} of $U_{r+1}(D)$ means the following.

\bclm{\bf Theorem.}\label{thmspeclz} {\em Let $A$ be a $\QQ$-algebra, $P\in A[T]$, $\phi=\sum_{n\geq
0}\phi_nt^n/n!\in A[[t]]$ and $(b_0,b_1,\ldots, b_r)\in A^{r+1}$ any $(r+1)$-tuple. Then  the unique
$\QQ$-algebra homomorphism, defined by $x_i\mapsto b_i$, $e_i\mapsto e_i(P)$ and $f_i\mapsto \phi_i$,
maps the universal solution $g$, as in~(\ref{eq:univsol}), to the unique solution of the Cauchy problem
\be \left\{\matrix{P(D)y&=&\phi,&\ \cr D^iy(0)&=&b_i\,,& 0\leq i\leq r.}\right.\label{eq:specialcauchy} \ee
 }
\eclm

For each $0\leq i\leq r$, let $\psi_i: E_r[\bfx,\bar\bff]\sra E_r$ be the unique $E_r$-algebra
homomorphism over the identity sending $\bfx\mapsto (\underbrace{0,\ldots, 0}_{i},1,h_1,\ldots,
h_{r-i})$ and $\bar\bff\mapsto (0,0,\ldots)$.
 \bclm{\bf Corollary.}\label{cor86} {\em If $u_i:=\psi_i(g)\in E_r[[t]]$, where
$g$ is the unique solution of  the universal Cauchy problem~(\ref{eq:cauchyprob}),  then
$\bfu=(u_0,u_1,\ldots, u_r) $ is an $E_r$-basis of $\ker U_{r+1}(D)$. } \eclm

\proof Using the same arguments as in Theorem~\ref{nhunsol}, one shows that $u_i$ is a solution of $U_{r+1}(D)y=0$. Furthermore, if $ u:=a_0u_0+a_1u_1+\ldots+a_ru_r=0 $, then $u$ is the unique solution of
$U_{r+1}(D)y=0$,  with the zero initial conditions. Then by uniqueness $u=0$,  i.e. $(u_0,\ldots, u_r)$
are linearly independent.\qed

\bclm{\bf Corollary.}\label{cor67} {\em Let $A$ be any $\QQ$-algebra and
$P\in A[T]$. Let $\psi: E_r\sra A$ be the unique morphism mapping $e_i\mapsto e_i(P)$. Then
$(\psi(u_0),\psi(u_1),\ldots, \psi(u_r))$ is an $A$-basis of $\ker P(D)$.\qed } \eclm

In other words, $\ker P(D)\cong \ker U_{r+1}(D)\otimes_{E_r} A$.

\claim{}  Let $n\geq 0$ be an integer and $\bmu$ a partition of length at most $r+1$  with weight $n$.
Denote by ${n\choose \bmu}$  the coefficient of $x_0^{\mu^0}x_1^{\mu^1}\ldots x_r^{\mu^r}$ in the
expansion of $(x_0+x_1+\ldots+x_r)^n$. With the usual convention  $0!=1$, one has
 \[
 {n\choose \bmu}={n!\over \mu_0!\mu_1!\ldots\mu_r!}.
 \]
In Section~\ref{intthgrbund} the Schur polynomials $\Delta_\bmu(a)$ associated to partition $\bmu$  and
to (the coefficients of) a formal power series $a=\sum_{n\geq 0} a_nt^n$ were defined, see
(\ref{eq:defschur}). In our notation, the coefficients form sequence $\bar\bfa=(a_0, a_1, \ldots\,)$;
below we will write $\Delta_\bmu(\bar\bfa)$ instead of $\Delta_\bmu(a)$.

 \bclm{\bf Theorem.}\label{rmk72} {\em For each partition $\blamb$, the following equality holds:
 \[
 W_\blamb(\bfu)=\sum_{n\geq 0}\sum_{|\bmu|=n}{n\choose \bmu}\Delta_{\blamb+\bmu}(\bar\bfh){t^n\over n!}.
 \]
 In particular, the "constant term" is $W_\blamb(\bfu)(0)=\Delta_\blamb(\bar\bfh)$. \qed}
 \eclm

 It is a straightforward combinatorial exercise made easy by the use of the basis $\bfu$ found in~\ref{cor86}.
 See~\cite{GatScherb1} for details.

 \bclm{\bf Proposition.}\label{giambthmwr} {\em Giambelli's formula for Wro\'nskians holds:
 \[
 W_\blamb(\bfu)=\Delta_\blamb(\bar\bfh)\cdot W_0(\bfu).
 \]
 }
 \eclm
\proof First of all, by Remark~\ref{rmkpropor},  $W_\blamb(\bfu)$ is proportional to $W_0(\bfu)$, i.e.
$W_\blamb(\bfu)=c_\blamb W_0(\bfu)$ for some $c_\blamb\in E_r$. Next, two formal power series are
proportional if and only if the coefficients of the same powers of $t$ are proportional, with the same
factor of proportionality. Finally,
\[
c_\blamb={W_{\blamb}(\bfu)(0)\over W_{0}(\bfu)(0)}=\Delta_\blamb({\bar\bfh}), \] according to
Theorem~\ref{rmk72}. \qed

\bclm{\bf Corollary.}\label{piericor} {\em Pieri's formula for generalized Wro\'nskians holds:
\begin{center}
$ h_iW_\blamb(\bfu)=\sum_\bmu W_\bmu(\bfu)$,
\end{center}
where the sum is over the partitions $\bmu=(\mu_0,\mu_1,\ldots, \mu_r)$  such that $|\bmu|=i+|\blamb|$
and
\[
\mu_0\geq\lambda_0\geq\mu_1\geq\lambda_1\geq\ldots\geq \mu_r\geq \lambda_r\,.
\]
}
\eclm
 It is well known that Giambelli's and Pieri's implies each other. See e.g.~\cite[Lemma A.9.4]{Fu1}.

\claim{} \label{fin813} Let now $\varrho_{r,d}: G\sra X$ be a Grassmann bundle, where $G:=G(r+1, F)$ and
$F$ is a vector bundle of rank $d+1$. As recalled in Section \ref{intthgrbund}, $A^*(G)$ is freely
generated as $A^*(X)$-module (see~\cite[Proposition 14.6.5]{Fu1}) by
\[
\Delta_\blamb(c_t(\Qcal_r-\varrho_{r,d}^*F))\cap[G].
\]
The exact sequence~(\ref{eq:univexseq}) implies that $c_t(\Scal_r)c_t(\Qcal_r)=c_t(\varrho_{r,d}^*F)$, which is
equivalent to
\[
1=c_t(\Scal_r)\,{c_t(\Qcal_r)\over c_t(\varrho_{r,d}^*F)}=c_t(\Scal_r)c_t(\Qcal_r-\varrho_{r,d}^*F).
\]
Set $\varepsilon_i=(-1)^{i}c_i(\Scal_r)$ and consider the differential equation \be
D^{r+1}y-\varepsilon_1\cdot D^ry+\ldots+(-1)^{r+1}\varepsilon_{r+1}\cdot y=0\,.\label{eq:schubdifeqII}
\ee We look for solutions in $(A^*(G)\otimes \QQ)[[t]]$. By Corollary~(\ref{cor67})  the unique morphism
$\psi : E_r\sra A^*(G)\otimes \QQ$, sending $e_i\mapsto \varepsilon_i$,  maps the universal fundamental
system $(u_0,u_1,\ldots, u_r)$ to $\bfv=(v_0,v_1,\ldots, v_r)$, where $v_i=\psi(u_i)$ and, as a
consequence, it maps $h_i$ to  $c_i(\Qcal_r-\varrho_{r,d}^*F)$ and $W_\blamb(\bfu)$ to $W_\blamb(\bfv)$.
Then we have proven that
\[
\Delta_\blamb(c_t(\Qcal_r-\varrho_{r,d}^*F))={W_\blamb(\bfv)\over W_0(\bfv)}\,.
\]
In other words,   the Chow group $A^*(G)$ can be identified with the $A^*(X)$-module generated by the
generalized Wro\'nskians associated to the basis $\bfv$ of solutions of the differential
equation~(\ref{eq:schubdifeqII}). In particular we have shown that the class
$[\Omega_{\blamb}(\varrho_{r,d}^*F_\bullet)]$ of the generalized Wro\'nski variety
$\Omega_{\blamb}(\varrho_{r,d}^*F_\bullet)$ is an $A^*(X)$-linear combination of ratios of generalized
Wro\'nskians associated to the basis $\bfv$ of ~(\ref{eq:schubdifeqII}), by virtue
of~(\ref{eq:gatsant}), ~(\ref{eq:relmuep}) and~(\ref{eq:kmpflksv}).

\end{document}